\documentclass{article}

\usepackage{quartics-layout}
\usepackage{quartics}

\usepackage[
	style=alphabetic,
	giveninits=true,
]{biblatex}
\title{The rational cohomology ring of the moduli space of GIT-stable plane quartic curves, and of the space of their equations}
\author{Maksym Dolgikh \\ Institut f\"ur Mathematik, Freie Universit\"at Berlin
	\and Alexey Gorinov \\ Faculty of Mathematics, HSE University}

\begin{document}

\maketitle

\begin{abstract}
	We study GIT-stable quartic curves in the complex projective plane, that is, plane quartics with at most nodes and cusps.
	We determine the rational cohomology rings of two spaces: the moduli space of such quartics and the space of their equations.
	The cohomology ring of the moduli space, which we compute together with its mixed Hodge structure, is generated by three classes that come from the classifying space of the general linear group, and we give an explicit ideal of relations.
	The relations are obtained from a single equivariant Euler class by integration along the fibre.
	In the cohomology ring of the space of equations, every non-trivial product involves the generator of degree five, and in particular this ring is not an exterior algebra, in contrast with the case of nodal quartics.
\end{abstract}

\section{Introduction}
\label{sec:intro}

A smooth quartic curve in $ \C P^2 $ is a canonically embedded non-hyperelliptic curve of genus three, so the moduli space of smooth plane quartics is the moduli space $ \mathcal{M}_3 $ of curves of genus three without its hyperelliptic locus.
Geometric invariant theory provides a natural partial compactification of this space, the moduli space of GIT-stable plane quartics, which is a geometric quotient.
A plane quartic is GIT-stable if and only if its only singular points are ordinary nodes and ordinary cusps.
This moduli space is the stable part of the GIT compactification of the moduli space of plane quartics, which is one of the log canonical models of the moduli space $ \overline{\mathcal{M}}_3 $ of stable curves of genus three \cite[Theorem 1]{HyeLee10}.

The rational cohomology of the moduli space of smooth plane quartics is known: it is $ \Q $ in degrees $ 0 $ and $ 6 $ and zero in all other degrees \cite{Loo93}, \cite[Example 2]{PetSte03}.
So is the cohomology ring of the moduli space of nodal plane quartics, which is $ \Q[a]/(a^4) $ with $ a $ of degree $ 2 $ \cite[Theorem 2]{BerGorKon23}.
In this paper we compute the rational cohomology ring of the moduli space of GIT-stable plane quartics, together with its mixed Hodge structure.
Our second result is the ring structure of the cohomology of the space of equations of such quartics.

Let $ \Pi = \C[x_0,x_1,x_2]_4 $ denote the $15$-dimensional complex vector space of homogeneous polynomials of degree $4$.
By \cite[Proposition 4.2]{mfk94}, a quartic in $\C P^2$ is \emph{not} GIT-stable if and only if it has a triple point or a tacnode.

\begin{definition}
	Let $ Q = \V(F) \subset \C P^2 $ be a curve corresponding to a polynomial $ F $, and let $ p \in \C P^2 $ be a point.
	\begin{itemize}[leftmargin=*]

		\item The point $ p \in \C P^2 $ is an \emph{$ m $-point} if in an affine chart centred at $ p $, the homogeneous components $ F_k $ of $ F $ of degree $ k $ are zero for all $ k < m $.

		\item A line $ l = \V(L) \subset \C P^2 $ is an \emph{$ m $-fold tangent at $ p $} if in an affine chart centred at $ p $, the homogeneous components $ F_k $ of $ F $ of degree $ k $ are multiples of $ L $ for all $ k < m $.

		\item A \emph{cusp} is a double point of a curve with only one $3$-fold tangent line, i.e. a point whose projective tangent cone consists of a single point.
		A \emph{(classical) tacnode} is a cusp whose $3$-fold tangent is also a $4$-fold tangent.
	\end{itemize}
\end{definition}

Let $\Sigma \subset \Pi$ be the discriminant consisting of polynomials corresponding to quartics with a triple point or a tacnode, and set
\begin{equation*}
	\UU \coloneq \Pi \setminus \Sigma .
\end{equation*}
Thus $ \UU $ is the space of equations of GIT-stable quartic curves.
The group $\GL_3(\C)$ acts on $\UU$ by linear changes of coordinates in $\C^3$.
All stabilisers of this action are finite, and the quotient $ \UU/\GL_3(\C) $ is a geometric quotient \cite[Appendix A]{BerGorKon23}.
It is the moduli space of GIT-stable plane quartics.

For a space $ X $, we denote by $ H^*(X) $ its rational singular cohomology, and by $ \overline H_*(X) $ its Borel--Moore homology.
Let $\Q(w)$ denote Tate mixed Hodge structure on $\Q$ of weight $-2w$.

Since the stabilisers are finite, the rational cohomology of $ \UU/\GL_3(\C) $ is the $ \GL_3(\C) $-equivariant cohomology of $ \UU $.
It therefore receives a ring homomorphism from $ H^*(\BGL_3(\C)) = \Q[s_1, s_2, s_3] $, where $ s_i $ has degree $ 2i $, and we denote the image of $ s_i $ again by $ s_i $.
Our main result is the following.

\begin{alphatheorem}[\autoref{thm:cohom_sp_of_moduli}, \autoref{thm:cohom_ring_moduli}]
	\label{thm:intro_b}
	The singular cohomology $ H^*( \UU/\GL_3(\C)) $ of the moduli space $  \UU/\GL_3(\C) $ can be described as follows.
	\begin{enumerate}[(1), leftmargin=*]
		\item As a graded vector space, together with its mixed Hodge structure, $ H^*( \UU/\GL_3(\C)) $ is given by the table below.
		\begin{equation*}
			{
				\begin{array}{c|ccccc}
					n                 & 0     & 2      & 4                 & 6                 & 8      \\ \hline
					H^n( \UU/\GL_3(\C)) & \Q(0) & \Q(-1) & \Q(-2)^{\oplus 2} & \Q(-3)^{\oplus 2} & \Q(-4)
				\end{array}
			}
		\end{equation*}
		\item As a ring, $ H^*( \UU/\GL_3(\C)) = \Q[s_1, s_2, s_3] / I $, where $ s_i \in H^{2i} $, and the relations $ I $ are given by
		\begin{gather*}
			0 = 72s_1^3 - 39s_1s_2 + 7s_3 \in H^6, \\
			0 = s_1^4 - \frac{25}{369}s_2^2 = s_1^2s_2 - \frac{17}{123}s_2^2 = s_1s_3-\frac{3}{41}s_2^2 \in H^8, \\
			0 = s_1s_2^2 \in H^{10} .
		\end{gather*}
	\end{enumerate}
\end{alphatheorem}

Part (2) says in particular that the homomorphism $ H^*(\BGL_3(\C)) \to H^*(\UU/\GL_3(\C)) $ is surjective.
In degrees up to $ 4 $ it is an isomorphism for a simple reason: the discriminant $ \Sigma $ has complex codimension three, and removing it does not change the equivariant cohomology in these degrees.
The theorem shows that the classes $ s_1, s_2, s_3 $ still generate the ring in the higher degrees, where the relations appear.
For smooth quartics this is not the case: the class of degree $ 6 $ in the cohomology of their moduli space has weight $ 12 $ \cite[Example 2]{PetSte03}, whereas a class of degree $ 6 $ that comes from $ \BGL_3(\C) $ has weight $ 6 $.

The relation of lowest degree has a geometric meaning.
All relations are obtained from the equivariant Euler class of one vector bundle by integration along the fibre, as explained in the strategy of the proofs below, and the image of this Euler class itself is $ -8(72s_1^3 - 39s_1s_2 + 7s_3) $.
A general element of $ \Sigma $ has a single tacnode, and it follows that this image is, up to sign, the $ \GL_3(\C) $-equivariant fundamental class of $ \Sigma $.
Accordingly, its restriction to the scalar matrices recovers the degree of the projectivisation of $ \Sigma $ in $ \PP^{14} $, which is $ 200 $, the classical number of plane quartics with a tacnode through eleven general points \cite[Table 1]{Zin05}.

Our second result concerns the space of equations $ \UU $ itself.
As a graded vector space with its mixed Hodge structure, $ H^*(\UU) $ is computed in \cite[Theorem A]{Dol1}, which is restated here as \autoref{thm:cohom_of_equations}.
Using \autoref{thm:intro_b}, we determine its ring structure.

\begin{alphatheorem}[\autoref{thm:eqs_mult}]
	\label{thm:intro_a}
	As a ring, $ H^*(\UU) $ is generated by classes
	\begin{equation*}
		x \in H^5, \quad y, y' \in H^7, \quad z \in H^9, \quad t \in H^{10}, \quad e \in H^{12} ,
	\end{equation*}
	with the only non-trivial products given by:
	\begin{equation*}
		xy = -yx, \quad xy'=-y'x , \quad xz=-zx , \quad xt=tx , \quad xe=ex .
	\end{equation*}
	Moreover, $ e, xy , xy' \in H^{12} $ are linearly independent.
\end{alphatheorem}

Thus every non-trivial product involves the class $ x $ of degree $ 5 $.
In particular, $ H^*(\UU) $ is not an exterior algebra, in contrast with the space of equations of nodal plane quartics, whose cohomology ring is an exterior algebra on generators of degrees $ 3 $, $ 5 $ and $ 7 $ \cite[Theorem 2]{BerGorKon23}.
Unlike the other five, the class $ e $ is not determined by \autoref{thm:intro_a}: the theorem holds for any $ e \in H^{12} $ with $ xe \neq 0 $ (\autoref{rem:e_normalisation}).
The table of \cite[Theorem A]{Dol1} and the products of \autoref{thm:intro_a} coincide with the Koszul homology of the ring of \autoref{thm:intro_b} over $ H^*(\BGL_3(\C)) $, with its bigrading and its products.
This observation is not used in the proofs, and it is recorded in \autoref{sec:koszul}.

This is the second of three papers.
The first \cite{Dol1} computes the graded vector space $ H^*(\UU) $ and its mixed Hodge structure by the method of conical resolutions, and this is the starting point of the present paper.
The third determines the cohomology ring of the space $ \UU/\C^{\times} $ of GIT-stable quartic curves, the morphisms induced by the natural quotient maps between $ \UU $, $ \UU/\C^{\times} $ and $ \UU/\GL_3(\C) $ and by the orbit maps, and the morphisms induced by the inclusion of the space of nodal quartics into the space of GIT-stable quartics.
The three papers are based on the first author's PhD thesis, written at Freie Universit\"at Berlin.
The present paper can be read on its own: what is used from \cite{Dol1} is restated where it is needed.

\subsection*{Related work}

The moduli space of GIT-stable plane quartics has several neighbours whose rational cohomology is known.
Looijenga computed the cohomology of $ \mathcal{M}_3 $ and of the moduli space of smooth plane quartics \cite{Loo93}.
Tommasi computed the cohomology of the moduli space of smooth curves of genus four \cite{Tom05} and of smooth curves of genus three with two marked points \cite{Tom07}, and Bergstr\"om and Tommasi computed that of the moduli space of stable curves of genus four \cite{BerTom07}.
For nodal plane quartics, the cohomology rings of the space of equations and of the moduli space were computed by Berdnikov, Gorinov and Konovalov \cite[Theorem 2]{BerGorKon23}.

Apart from Looijenga's, these computations use the Vassiliev method of conical resolutions \cite{Vass99, Gor05}, and so does \cite{Dol1}, on which the present paper rests.
Further applications of the method are listed in the introduction of \cite{Dol1}.

The moduli space of GIT-stable plane quartics is an open subset of the GIT compactification of the moduli space of plane quartics.
In a recent preprint, Ascher and Lee study the Chow and cohomology rings of the GIT, K-moduli and Hacking moduli stacks of plane quartics, and compute in particular the intersection Betti numbers of the GIT compactification \cite[Theorem 7.6]{AscLee26}.

The technique of this paper has a counterpart in intersection theory, where equivariant excision of a discriminant is used to compute Chow rings of stacks of plane curves: of smooth plane quartics \cite{DiLFulVis21}, of nodal plane curves of every degree \cite{CelUruYan25}, and of plane quartics with singularities of type $ A $ \cite[Corollary 3.15]{Per23c}.
In \autoref{sec:chow} we show that the cycle class map from the rational Chow ring of the quotient stack $ [\UU/\GL_3(\C)] $ to $ H^*(\UU/\GL_3(\C)) $ is an isomorphism, so that \autoref{thm:intro_b} also presents this Chow ring.
A presentation of the rational Chow ring of the stack of GIT-stable plane quartics was obtained independently by Ascher and Lee \cite[Theorem 9.15]{AscLee26}, and it agrees with ours, as we check at the end of \autoref{sec:chow}.

\subsection*{Strategy of the proofs}

We prove part $(1)$ of \autoref{thm:intro_b} in \autoref{sec:equivariant_quotient}, where we determine the additive structure of $ H^*(\UU/\GL_3(\C)) $ by comparing two spectral sequences.

To prove part $(2)$ of \autoref{thm:intro_b}, we show that $H^*(\UU/ \GL_3(\C)) =  H^*( \UU_{\hGL_3(\C)} )$ is isomorphic to the ring \mbox{$  H^*(\BGL_3(\C)) = \Q[s_1, s_2, s_3] $} modulo the image under integration along the fibre map of an ideal generated by the equivariant Euler class of a certain fibre bundle (see \autoref{sec:simplifying_to_int_along_fibre}).
We use the Leray--Hirsch theorem together with \autoref{prop:int_fibre_explicitly} to make the integration along the fibre map explicit in terms of $ s_i $, and compute the relations on $ s_i $ using SageMath \cite{sage}.

This, in turn, allows us to establish \autoref{thm:intro_a} using the multiplicative structure on the Leray spectral sequence of the quotient map $ \UU \to \UU / \GL_3(\C) $, which is studied in \autoref{sec:mult_struct_eqs} and whose limiting term can be found in \autoref{table:GE_inf_gens}.

\section{Additive structure of \texorpdfstring{$  H^*(\UU_{\hGL_3(\C)};\Q) $}{}}
\label{sec:equivariant_quotient}

For brevity, we write $ G \coloneq \GL_3(\C) $.

Our starting point is the cohomology of the space of equations $ \UU $.
\begin{theorem}[{\cite[Theorem~A]{Dol1}}]\label{thm:cohom_of_equations}
	As a graded vector space, together with its mixed Hodge structure, the singular cohomology $ H^*(\UU) $ of the space $ \UU $ of equations of GIT-stable quartics is given by the table below.
	\begin{equation*}
		\resizebox{0.95\textwidth}{!}{$
			\begin{array}{c|ccccccccc}
				n      & 0     & 5      & 7                 & 9      & 10     & 12                & 14     & 15     & 17      \\ \hline
				H^n(\UU) & \Q(0) & \Q(-3) & \Q(-4)^{\oplus 2} & \Q(-5) & \Q(-6) & \Q(-7)^{\oplus 3} & \Q(-8) & \Q(-9) & \Q(-10)
			\end{array}
		$}
	\end{equation*}
\end{theorem}

We use it to compute the cohomology $ H^*(\UU / G) $ of the geometric quotient and $ H^*(\UU_{hG}) $ of the equivariant quotient.
Furthermore, in our case $ H^*(\UU / G) = H^*(\UU_{hG}) $ (for details, see \cite[Appendix A]{BerGorKon23}).
We begin by determining the additive structure of $ H^*(\UU / G) $ by comparing two spectral sequences.
Determining the multiplicative structure requires considerably more work.
The key idea is to first reduce the problem to finding the images of differentials in a Borel--Moore homological spectral sequence, and then to actually calculate said images.
This final step is carried out using SageMath \cite{sage}, a computer algebra system.
Using the obtained ring structure on $ H^*(\UU / G) $, we will later find the ring structure on $ H^*(\UU) $.

The vector spaces $ H^j(\UU_{hG}) = H^j(\UU/G) $ can be determined by comparing the cohomological spectral sequence $ {}^{\UU}E^{}_{*} $ of the fibre bundle
\begin{equation*}
	\UU \to \UU_{hG} \to BG ,
\end{equation*}
and the cohomological spectral sequence $ {}^{G}E^{}_{*} $ corresponding to
\begin{equation*}
	G \to \UU \to \UU/G .
\end{equation*}
Although the latter need not be a fibre bundle, the associated Leray spectral sequence behaves as if it were.
In particular, by \cite[Corollary A.1]{BerGorKon23}, the second page is a tensor product of its $ \nth{0} $ row and $ \nth{0} $ column.

\begin{theorem}\label{thm:cohom_sp_of_moduli}
	As a graded vector space, together with its mixed Hodge structure, $ H^*( \UU/\GL_3(\C)) $ is given by the table below.
	\begin{equation*}
		{
			\renewcommand{\arraystretch}{1.4}
			\setlength{\tabcolsep}{8pt}
			\begin{array}{c|ccccc}
				n                 & 0     & 2      & 4                 & 6                 & 8      \\ \hline
				H^n( \UU/\GL_3(\C)) & \Q(0) & \Q(-1) & \Q(-2)^{\oplus 2} & \Q(-3)^{\oplus 2} & \Q(-4)
			\end{array}
		}
	\end{equation*}
\end{theorem}

\begin{proof}
	By \cite[Section 22.1]{fomenko2016homotopical}, $ {}^{\UU}E_2^{pq} = H^p(BG; H^q( \UU )) $.
	For $ G = \GL_3(\C) $, we have \mbox{$ H^*(G; \Q) = \Lambda [
	\overset{1}{\xi_1},
	\overset{3}{\xi_3},
	\overset{5}{\xi_5}
	] $}, where $ \xi_i \in H^{i}(G; \Q) $.
	Recall also that $ H^*(BG; \Q) = \Q[s_1, s_2, s_3] $, where $ s_i \in H^{2i}(BG; \Q) $.
	See \cite[\textsection 1]{BGG73} for details.

	Thus the bottom left corner of $ {}^{\UU}E^{}_{2} $ is as in the table \autoref{table:UE2}.

	\begin{table}[!htbp]
		\centering
		\caption{Bottom left corner of ${}^{\UU}E^{}_{2}$}
		\label{table:UE2}
		\begin{tabular}{|r|c|c|c|c|c|c|c|c|c|}\hline
			& 0                                   & 1 & 2                                  & 3 & 4           & 5 & 6                                   & 7 & 8                                   \\
			\hline
			7 &  $ 2\Q(-4) $ &   & $ 2\Q(-5) $                        &   & $ 4\Q(-6) $ &   & $ 6\Q(-7) $                         &   & $ 8\Q(-8) $                         \\
			\hline
			6 &                                     &   &                                    &   &             &   &                                     &   &                                     \\
			\hline
			5 &  $ \Q(-3) $  &   &  $ \Q(-4) $ &   & $ 2\Q(-5) $ &   & $ 3\Q(-6) $                         &   & $ 4\Q(-7) $                         \\
			\hline
			4 &                                     &   &                                    &   &             &   &                                     &   &                                     \\
			\hline
			3 &                                     &   &                                    &   &             &   &                                     &   &                                     \\
			\hline
			2 &                                     &   &                                    &   &             &   &                                     &   &                                     \\
			\hline
			1 &                                     &   &                                    &   &             &   &                                     &   &                                     \\
			\hline
			0 & $ \Q(0) $                           &   & $ \Q(-1) $                         &   & $ 2\Q(-2) $ &   &  $ 3\Q(-3) $ &   &  $ 4\Q(-4) $ \\
			\hline
			& 0                                   & 1 & 2                                  & 3 & 4           & 5 & 6                                   & 7 & 8                                   \\
			\hline
		\end{tabular}
	\end{table}

	Note that $ \dim \UU/G = \dim \UU - \dim G = 2\cdot 15 - 2\cdot 9 = 12 $, so $ H^j(\UU_{hG}) = H^j(\UU/G) = 0 $ for $ j > 12 $.
	Moreover, $ H^{12}(\UU_{hG}) = 0 $.
	Indeed, consider the corner $ (p,q) = (12,9) $ of \mbox{$ {}^{G}E_2^{pq} = H^p(\UU/G) \otimes H^q(G) $}.
	For every $ r \geq 2 $, the differential into that spot starts at $ (12-r, 8+r) $, in a row $ \geq 10 $, where $ H^{q}(G) = 0 $ because $ \Lambda[\xi_1,\xi_3,\xi_5] $ vanishes above degree $ 9 $.
	And the differential out of it lands at $ (12+r, 10-r) $, in a column $ \geq 14 $, where $ H^{p}(\UU/G) = 0 $ because $ \dim \UU/G = 12 $.
	Hence $ {}^{G}E^{12,9}_{\infty} = {}^{G}E^{12,9}_{2} = H^{12}(\UU/G) \otimes \langle \xi_1\xi_3\xi_5 \rangle $, which is a subquotient of $ H^{21}(\UU) = 0 $ (see \autoref{thm:cohom_of_equations}), so $ H^{12}(\UU/G) = 0 $.
	Therefore, the second page of $ {}^{G}E^{}_{*} $ is as in \autoref{table:partial_GE2}.

	\begin{table}[!htbp]
		\centering
		\caption{Known part of ${}^{G}E^{}_{2}$}
		\label{table:partial_GE2}
		\begin{tabular}{|r|C|C|C|C|C|C|C|C|C|C|C|C|}\hline
			& 0      & 1 & 2      & 3 & 4        & 5 & 6 & 7 & 8 & 9 & 10 & 11 \\ \hline
			9 & \Q(-6) &   & \Q(-7) &   & 2 \Q(-8) & ? & ? & ? & ? & ? & ?  & ?  \\ \hline
			8 & \Q(-5) &   & \Q(-6) &   & 2 \Q(-7) & ? & ? & ? & ? & ? & ?  & ?  \\ \hline
			7 &        &   &        &   &          &   &   &   &   &   &    &    \\ \hline
			6 & \Q(-4) &   & \Q(-5) &   & 2 \Q(-6) & ? & ? & ? & ? & ? & ?  & ?  \\ \hline
			5 & \Q(-3) &   & \Q(-4) &   & 2 \Q(-5) & ? & ? & ? & ? & ? & ?  & ?  \\ \hline
			4 & \Q(-3) &   & \Q(-4) &   & 2 \Q(-5) & ? & ? & ? & ? & ? & ?  & ?  \\ \hline
			3 & \Q(-2) &   & \Q(-3) &   & 2 \Q(-4) & ? & ? & ? & ? & ? & ?  & ?  \\ \hline
			2 &        &   &        &   &          &   &   &   &   &   &    &    \\ \hline
			1 & \Q(-1) &   & \Q(-2) &   & 2 \Q(-3) & ? & ? & ? & ? & ? & ?  & ?  \\ \hline
			0 & \Q(0)  &   & \Q(-1) &   & 2 \Q(-2) & ? & ? & ? & ? & ? & ?  & ?  \\ \hline
			& 0      & 1 & 2      & 3 & 4        & 5 & 6 & 7 & 8 & 9 & 10 & 11 \\ \hline
		\end{tabular}
	\end{table}

	The spectral sequence converges to $ H^*(\UU) $.
	Since the group $ H^{i}(\UU) = 0 $ is zero for each \mbox{$ i > 17 $}, the groups \mbox{$ {}^{G}E^{10,9}_{2} = 0 = {}^{G}E^{11,9}_{*} $} must be $ 0 $, as no differential can kill them.
	Consequently, \mbox{$ {}^{G}E^{10,q}_{2} = 0 = {}^{G}E^{11,q}_{*} $} are zero for all $ q $.
	In particular, we obtain that the groups \mbox{$ H^{10}(\UU_{hG}) = H^{10}(\UU/G) = 0 = H^{11}(\UU/G) = H^{11}(\UU_{hG}) $} in the $\nth{0}$ row are zero.

	But then for the same reason ${}^{G}E^{9,q}_{2} = 0$ for all $ q $, and in particular $ H^{9}(\UU_{hG})=0 $.

	As $ H^{17}(\UU) = \Q(-10) $, we have $ {}^{G}E^{8,9}_{2} = \Q(-10)$, and therefore $ H^{8}(\UU_{hG}) = \Q(-4) $.

	As $ H^{16}(\UU) = 0 $, ${}^{G}E^{7,q}_{2} = 0$, and in particular $ H^{7}(\UU_{hG})=0 $.
	Moreover, the differential $ \dif^{6,9}_2 \colon {}^{G}E^{6,9}_{2} \to {}^{G}E^{8,8}_{2} $ must be surjective.
	As $ H^{15}(\UU) = \Q(-9) $, we obtain $ \ker \dif^{6,9}_2 = \Q(-9) $.
	It follows that $ {}^{G}E^{6,9}_{2} = 2\Q(-9) $, and therefore $ H^{6}(\UU_{hG}) = 2\Q(-3)$.

	Finally, $ H^{13}(\UU) = 0 $, but $ {}^{G}E^{5,8}_{*} $ cannot be killed by any differential.
	Hence $ {}^{G}E^{5,8}_{2} = 0 $, and therefore $ H^{5}(\UU_{hG}) = 0 $.

	Now we know the additive structure of $ H^{*}(\UU_{hG}) $, and we can see that $ {}^{G}E^{}_{2} $ is as in \autoref{table:GE_2}
	\begin{table}[!htbp]
		\centering
		\caption{Full ${}^{G}E^{}_{2}$}
		\label{table:GE_2}
		\begin{tabular}{|r|C|C|C|C|C|C|C|C|C|}
			\hline
			& 0      & 1 & 2      & 3 & 4        & 5 & 6       & 7 & 8       \\ \hline
			9 & \Q(-6) &   & \Q(-7) &   & 2 \Q(-8) &   & 2\Q(-9) &   & \Q(-10) \\ \hline
			8 & \Q(-5) &   & \Q(-6) &   & 2 \Q(-7) &   & 2\Q(-8) &   & \Q(-9)  \\ \hline
			7 &        &   &        &   &          &   &         &   &         \\ \hline
			6 & \Q(-4) &   & \Q(-5) &   & 2 \Q(-6) &   & 2\Q(-7) &   & \Q(-8)  \\ \hline
			5 & \Q(-3) &   & \Q(-4) &   & 2 \Q(-5) &   & 2\Q(-6) &   & \Q(-7)  \\ \hline
			4 & \Q(-3) &   & \Q(-4) &   & 2 \Q(-5) &   & 2\Q(-6) &   & \Q(-7)  \\ \hline
			3 & \Q(-2) &   & \Q(-3) &   & 2 \Q(-4) &   & 2\Q(-5) &   & \Q(-6)  \\ \hline
			2 &        &   &        &   &          &   &         &   &         \\ \hline
			1 & \Q(-1) &   & \Q(-2) &   & 2 \Q(-3) &   & 2\Q(-4) &   & \Q(-5)  \\ \hline
			0 & \Q(0)  &   & \Q(-1) &   & 2 \Q(-2) &   & 2\Q(-3) &   & \Q(-4)  \\ \hline
			& 0      & 1 & 2      & 3 & 4        & 5 & 6       & 7 & 8       \\ \hline
		\end{tabular}
	\end{table}
\end{proof}

To determine the multiplicative structure of $ H^{*}(\UU_{hG}) $, further work is required.

\section{Multiplicative structure of \texorpdfstring{$ H^{*}(\UU_{\hGL_3(\C)};\Q) $}{}}
\label{sec:v1_mult_equivariant}

Consider the fibre bundle
\begin{equation*}
	\UU  \to \UU_{hG} \to BG .
\end{equation*}
Recall that $ {}^{\UU}E_r $ is the corresponding cohomology spectral sequence.

By looking at the graded vector space $ H^{*}(\UU_{hG}) $ that we just calculated, we see that $ {}^{\UU}E^{}_{\infty} $ has exactly one non-zero row: $\nth{0}$.
Thus we can use the generators of $ H^*(BG; \Q) = {}^{\UU}E^{*0}_{2} $ and relations on them coming from the differentials to find the multiplicative structure on the limit $ H^{*}(\UU_{hG}) $.
For this, we need to know the images of the differentials that map to $\nth{0}$ row explicitly.

We proceed as follows.
First, we will show that $ {}^{\UU}E^{}_{*} $ is isomorphic to a BM-homological spectral sequence $ {}^{\UU}E^{*}_{} $.
Next, we will show how this isomorphism helps us describe some of the images of the differentials in the original spectral sequence as the image of a spectral sequence morphism.
We will later see that, for dimension reasons, \emph{all} of the images of the differentials can be described this way.
We will then show how to describe those images in terms of integration along the fibre and the equivariant Euler class of a vector bundle.
Finally, we compute them explicitly in terms of polynomials.

\subsection{A Borel--Moore homology spectral sequence}\label{sec:supp_ss_for_mult}

Let us establish some notation.
As before, we write $ G = \GL_3(\C) $ for brevity.
We want to work with $ \overline{H}^G_*(Y) = \overline{H}_*(Y_{hG}) $, but it might not make sense literally, as $ Y_{hG} $ may be not compact (finitely dimensional).
This is solved by replacing $ EG $ with a finite-dimensional model.
That is, a space whose low dimension cohomology groups coincide with those of $ EG $.
Equivalently, their high dimension Borel--Moore homology groups coincide.
The rest could be different.

Let $ BG^{(i)} \coloneq \Gr(3,\C^{i}) $ be the Grassmannian of $3$-dimensional linear subspaces in the complex vector space $ \C^{i} $.
Then $ BG = \varinjlim BG^{(i)} $ is the colimit of $ BG^{(i)} $, and the natural morphisms \mbox{$ \pi_{\leq j(i)}(BG^{(i)}) \to \pi_{\leq j(i)}(BG) $} are isomorphism, where $ i $ can be chosen to make $ j(i) $ arbitrarily large.
Let $ EG^{(i)} $ and $ Y^{(i)}_{hG} $ be the preimages of $ BG^{(i)} $ and $ Y_{hG} $ under the corresponding projection maps (here, $ Y $ could be $ \UU $ or $ \Pi $, for example).
Then $ EG = \varinjlim EG^{(i)} $ (with the topology of the direct limit) and $ Y_{hG} = \varinjlim Y^{(i)}_{hG} $.

By choosing large enough $ i_0 $, we can get isomorphisms
\[ H^j(BG^{(i_0)}) = H^j(BG) , \qquad H^j(\UU^{(i_0)}_{hG}) = H^j(\UU_{hG}) \]
for all $ j < j_0 $, with $ j_0 $ as large as we need.
So for large enough $ i_0 $, the spectral sequence of the fibre bundle
\begin{equation*}
	\UU  \to \UU^{(i_0)}_{hG} \to BG^{(i_0)} .
\end{equation*}
has the exact same bottom left corner as our spectral sequence $ {}^{\UU}E_*^{} $, and by the remarks on dimension from the previous section, it converges to the cohomology groups of our space $ \UU_{hG} $ up to degree $ j_0 $.

Let us fix a good enough $ i_0 $, and let us denote by $ {}^{\UU}E_*^{} $ the cohomological spectral sequence of the ``finite model'' fibre bundle.
From now on, we write $ \widetilde{Y} $ for the finite model $ Y^{(i_0)} $ of a space $ Y $.
Finally, let us denote $ \overline{H}^{G}_{*}(-) \coloneq \overline{H}^{}_{*}(- \times_G \widetilde{EG}) $, where $ \widetilde{EG} = EG^{(i_0)} $ is, as before, a ``good enough'' fixed model of $ EG $.
This way it makes sense to talk about $ \overline{H}^{G}_{*}(Y) = \overline{H}_{*}(\widetilde{Y_{hG}}) $.

Let $ {}^{\UU}E^* $ be the BM homology spectral sequence corresponding to the ``finite model'' fibre bundle above.
As we will soon see, $ {}^{\UU}E^{}_{*} $ is isomorphic to $ {}^{\UU}E^{*}_{} $, and $ {}^{\UU}E^{*}_{} $ is easier to work with.
Note that $ {}^{\UU}E^2_{pq} = \overline{H}_p(\widetilde{BG}; \overline{H}_q( \UU )) $.

We are now ready to build an isomorphism between spectral sequences $ {}^{\UU}E_* $ and $ {}^{\UU}E^* $.

\begin{remark}
	A cohomological spectral sequence $ E_* $ is a homological spectral sequence via \mbox{$ E_r^{pq} = E^r_{-p,-q} $}.
\end{remark}

Let $ N \coloneq \dim \widetilde{BG} $, $ d \coloneq \dim( \UU  ) = 2 \dim_{\C}(\Pi ) = 30 $.
Now, $ {}^{\UU,t}{E}_r^{pq} \coloneq {}^{\UU}E^r_{N - p, d - q} $ is a cohomological spectral sequence.

\begin{proposition}
	For a sufficiently large finite-dimensional approximation $ \widetilde{BG} $, the Borel--Moore homology spectral sequence
	$
	{}^{\UU}E^2_{pq} = \overline{H}_p(\widetilde{BG}; \overline{H}_q(\UU))
	$
	is isomorphic, from the $E^2$-page onwards, to the cohomological spectral sequence associated with the fibre bundle
	$
	\UU \to \widetilde{\UU_{hG}} \to \widetilde{BG}.
	$
	Moreover, after reindexing via
	$
	{}^{\UU,t}E_r^{pq} = {}^{\UU}E^r_{N-p,\, d-q},
	$
	this isomorphism is induced by Poincar\'e duality on $\UU$ and $\widetilde{BG}$.
\end{proposition}

\begin{proof}
	By Poincar\'e duality, $ PD_{\UU} \colon H^*( \UU ) \to \overline{H}_{d -*}( \UU ) $ and $ PD_{\widetilde{BG}} \colon H^*( \widetilde{BG} ) \to \overline{H}_{N -*}( \widetilde{BG} ) $ are isomorphisms.
	We show that they give an isomorphism of cohomological spectral sequences \mbox{$ PD \colon {}^{\UU}E_r \to {{}^{\UU,t}{E}}_r $}.
	Write $ \pi \colon \widetilde{\UU_{hG}} \to \widetilde{BG} $ for the projection, fix a smooth triangulation $ K $ of $ \widetilde{BG} $, and let $ K' $ be its barycentric subdivision and $ K^* $ its dual cell decomposition, whose $ p $-skeleton $ K^{*(p)} $ is the union of the cells dual to the simplices of dimension at least $ N-p $.
	For each of these cell structures $ L $, let $ {}^{\UU}E_*(L) $ be the spectral sequence in singular cohomology of the filtration of $ \widetilde{\UU_{hG}} $ by the preimages $ \pi^{-1}(L^{(p)}) $ of the skeleta.
	Then $ {}^{\UU}E_* = {}^{\UU}E_*(K) $, and $ {}^{\UU}E^* $ is the spectral sequence of the same filtration in Borel--Moore homology.

	We use these cell structures as follows.
	Poincar\'e duality on $ \widetilde{\UU_{hG}} $ turns the filtration by the preimages of the skeleta of $ K $ in Borel--Moore homology into a filtration by open sets in cohomology, and these open sets are homotopy equivalent to the preimages of the skeleta of $ K^* $.
	As $ K' $ refines both $ K $ and $ K^* $, it then lets us compare $ {}^{\UU}E_*(K^*) $ with $ {}^{\UU}E_*(K) $.

	First, we identify $ {}^{\UU,t}E_* $ with a spectral sequence in cohomology.
	Poincar\'e duality on an oriented $ n $-manifold $ X $ is the cap product $ \alpha \mapsto \alpha \frown [X] $ with the fundamental class $ [X] \in \overline{H}_n(X) $, and for an open subset $ V \subseteq X $ the class $ [X] $ restricts to $ [V] $.
	Capping therefore commutes with restriction to $ V $ in cohomology and with the corresponding map $ \overline{H}_*(X) \to \overline{H}_*(V) $ in Borel--Moore homology.
	Apply this to the open sets \mbox{$ W_p \coloneq \widetilde{\UU_{hG}} \setminus \pi^{-1}(K^{(N-1-p)}) $}, and use the Borel--Moore exact sequence of a closed subspace and its open complement.
	Capping with $ [W_p] $ gives isomorphisms $ H^m(W_p, W_{p'}) \to \overline{H}_{N+d-m}(W_p \setminus W_{p'}) $ for $ p' \leq p $ which, after a change of sign depending only on $ m $, commute with the maps of the two long exact sequences.
	So $ {}^{\UU,t}E_* $ is isomorphic, from $ E_1 $ on, to the spectral sequence in cohomology of the filtration of $ \widetilde{\UU_{hG}} $ by the open sets $ W_p $.

	Next, we show that this spectral sequence is isomorphic to $ {}^{\UU}E_*(K^*) $.
	In $ K' $, the complement of $ K^{(N-1-p)} $ deformation retracts onto the full subcomplex spanned by the barycentres of the simplices of dimension at least $ N-p $ (push the other barycentric coordinates linearly to zero), and this subcomplex is $ K^{*(p)} $.
	As $ \pi $ is a fibration, the inclusions $ \pi^{-1}(K^{*(p)}) \subseteq W_p $ are homotopy equivalences, so by the five lemma restriction is an isomorphism from the spectral sequence of the $ W_p $ to $ {}^{\UU}E_*(K^*) $.

	Finally, we compare $ {}^{\UU}E_*(K^*) $ with $ {}^{\UU}E_*(K) $ through $ K' $.
	Both $ K^{(p)} $ and $ K^{*(p)} $ lie in $ K'^{(p)} $, so restriction gives morphisms $ {}^{\UU}E_*(K) \leftarrow {}^{\UU}E_*(K') \to {}^{\UU}E_*(K^*) $, which on $ E_2 $ are the maps of cellular cohomology $ H^p(-; H^q(\UU)) $ induced by the identity of $ \widetilde{BG} $, hence isomorphisms.
	Composing, and as a morphism of spectral sequences that is an isomorphism on $ E_2 $ is one on every later page, we obtain $ PD \colon {}^{\UU}E_r \to {}^{\UU,t}E_r $ for $ r \geq 2 $.

	It remains to show that on $ E_2 $ the isomorphism $ PD $ is induced by $ PD_{\widetilde{BG}} $ and $ PD_{\UU} $.
	On $ E_1 $, each of the maps composed above is computed cell by cell, over contractible sets where the bundle is trivial.
	Under the identifications of the $ E_1 $ terms with cellular cochains or chains of $ \widetilde{BG} $ with coefficients $ H^q(\UU) $ or $ \overline{H}_{d-q}(\UU) $, it is therefore the corresponding map for the bundle $ \widetilde{BG} \to \widetilde{BG} $ with fibre a point, tensored with the identity, or with $ PD_{\UU} $ for the cap product up to a sign depending only on $ (p,q) $.
	For fibre a point, all these spectral sequences have a single non-zero row, so from $ E_2 $ on they are their abutments $ H^*(\widetilde{BG}) $ and $ \overline{H}_*(\widetilde{BG}) $, and the composite is the map of abutments, which is the cap product with $ [\widetilde{BG}] $ up to sign.
	Hence on $ E_2 $ the isomorphism $ PD $ is $ \pm PD_{\widetilde{BG}} \otimes PD_{\UU} $, the sign depending only on $ (p,q) $.
	Here $ \UU $, $ \widetilde{BG} $ and $ \widetilde{\UU_{hG}} $ are all complex manifolds, hence canonically and compatibly oriented.

	\begin{equation*}
		\begin{tikzcd}
			{}^{\UU}E^{pq}_{2} \arrow[r] \arrow[d]           & {{}^{\UU,t}E^{pq}_{2}} \arrow[d]                                   \\
			H^p(\widetilde{BG}) \otimes H^q( \UU ) \arrow[r] & \overline H_{N-p}(\widetilde{BG}) \otimes \overline H_{d-q}( \UU )
		\end{tikzcd}
	\end{equation*}
\end{proof}

The upper right corner of $ {}^{\UU}E^2 $, corresponding to the bottom left corner of $ {}^{\UU}E_2 $ that we've seen earlier, is shown in \autoref{table:UE2_bm_homol}.

\begin{table}[!htbp]
	\centering
	\caption{Upper right corner of $ {}^{\UU}E^2 $}
	\label{table:UE2_bm_homol}
	\begin{tabular}{|c|c|c|c|c|c|c|c|c|c|}
		\hline
		      & $N-8$       & $N-7$ & $N-6$       & $N-5$ & $N-4$       & $N-3$ & $N-2$       & $N-1$ & $N$         \\ \hline
		$d$   & $ 4\Q $ &       & $ 3\Q $ &       & $2\Q $  &       & $ \Q $  &       & $ \Q $   \\ \hline
		$d-1$ &             &       &             &       &             &       &             &       &             \\ \hline
		$d-2$ &             &       &             &       &             &       &             &       &             \\ \hline
		$d-3$ &             &       &             &       &             &       &             &       &             \\ \hline
		$d-4$ &             &       &             &       &             &       &             &       &             \\ \hline
		$d-5$ & $ 4\Q $ &       & $ 3\Q $ &       & $ 2\Q $ &       & $ \Q $  &       & $ \Q $  \\ \hline
		$d-6$ &             &       &             &       &             &       &             &       &             \\ \hline
		$d-7$ & $ 8\Q $ &       & $ 6\Q $ &       & $ 4\Q $ &       & $ 2\Q $ &       & $ 2\Q $ \\ \hline
			  & $N-8$       & $N-7$ & $N-6$       & $N-5$ & $N-4$       & $N-3$ & $N-2$       & $N-1$ & $N$         \\ \hline
	\end{tabular}
\end{table}

\subsection{Reduction to computing a morphism of spectral sequences}

We now apply the spectral sequence constructed above as follows.
We want to find the images of differentials mapping to the $\nth{0}$ row of $ {}^{\UU}E_* $ explicitly, so that we can deduce the multiplicative structure on the limit.
As shown in the previous section, this is equivalent to determining the images of differentials mapping to the $ d^{\text{th}} $ row of $ {}^{\UU}E^* $.
To determine these images, we show that there exists another spectral sequence $ {}^{\Pi}E^* $ such that the classes in $ {}^{\UU}E^2_{p,d} $ that vanish in the limit lie in the image of $ {}^{\Pi}E^2_{p,d} $ under a morphism of spectral sequences.
We will describe some of the classes in this image explicitly.
Subsequently, using the known additive structure, we show that no other classes in $ {}^{\UU}E^*_{p,d} $ vanish.

We will make use of the following lemma.
\begin{lemma}\label{lem:open_inclusion_spse}
	Let $ f \colon E_a \to E_b $ be an open inclusion of locally trivial fibre bundles over a common base space $ B $, and let $ \pi_a \colon E_a \to B $ and $ \pi_b \colon E_b \to B $ denote the projection maps.

	Then $f$ induces a morphism of the associated Borel--Moore homology spectral sequences \mbox{$\overline E^{*}(\pi_b) \longrightarrow \overline E^{*}(\pi_a)$},	which is compatible with the induced map $ f^{!} \colon \overline H_{*}(E_b) \to \overline H_{*}(E_a) $ on the abutments.

	In particular, $f^{!}$ preserves the filtrations on $\overline H_{*}(E_b)$ and $\overline H_{*}(E_a)$ arising from the spectral sequences, and the induced map
	\[
	F_p \overline H_{p+q}(E_b) / F_{p-1} \overline H_{p+q}(E_b) \longrightarrow F_p \overline H_{p+q}(E_a) / F_{p-1} \overline H_{p+q}(E_a)
	\]
	coincides with the map on $E^\infty$-pages determined by the morphism of spectral sequences.
\end{lemma}

The main result of this section is the following proposition.

\begin{proposition}
	Let $ {}^{\Pi}E^* $ denote the Borel--Moore homology spectral sequence associated with the fibre bundle
	$
		\Pi \to \widetilde{\Pi_{hG}} \to \widetilde{BG},
	$
	where $ \Pi = \C[x_0,x_1,x_2]_4 $.
	Let
	$
	f^2_{p,d} \colon {}^{\Pi}E^2_{p,d} \to {}^{\UU}E^2_{p,d}
	$
	be the morphism induced by the inclusion $ \UU \hookrightarrow \Pi $, and let $ f^{!} \colon \overline{H}^{G}_{*}(\Pi) \to \overline{H}^{G}_{*}( \UU ) $ be the induced map on the abutments.

	The sequence $ {}^{\Pi}E^* $ degenerates at the $ E^2 $-page, so that $ {}^{\Pi}E^2_{p,d} = \overline{H}^{G}_{p+d}(\Pi) $.
	Under this identification, every class in $ f^2_{p,d}( \ker f^{!} ) \subseteq {}^{\UU}E^2_{p,d} $ vanishes in the limit.
\end{proposition}

By exactness, $ \ker f^{!} $ is the image of $ g_* \colon \overline{H}^{G}_{*}(\Sigma) \to \overline{H}^{G}_{*}(\Pi) $, so the proposition applies to exactly the classes $ f^2_{p,d}(x) $ with $ x = g_*(y) $.
Note that the hypothesis $ x \in \ker f^{!} $ cannot be dropped: in row $ q = d $ the morphism $ f^2_{p,d} $ is induced on coefficients by the restriction $ \overline{H}_d(\Pi) \to \overline{H}_d( \UU ) $, which carries the fundamental class of $ \Pi $ to that of its dense open subset $ \UU $ and is therefore an isomorphism, so the image of $ f^2_{p,d} $ is all of $ {}^{\UU}E^2_{p,d} $.

Every class in ${}^{\UU}E^2_{p,d}$ that vanishes in the spectral sequence arises in this way, as will be shown later.

\begin{proof}
	By \cite[Section 22.1]{fomenko2016homotopical}, the spectral sequence $ {}^{\Pi}E^* $ satisfies
	\[
	{}^{\Pi}E^2_{pq} = \overline{H}_p(\widetilde{BG}; \overline{H}_q(\Pi)).
	\]
	Since $ \Pi = \C[x_0,x_1,x_2]_4 $ is a complex vector space of complex dimension $15$, we have \mbox{$ d = 2 \dim_{\C}(\Pi) = 30 $}, and
	\begin{equation*}
		\overline{H}_q(\Pi ; R) =
		\begin{cases}
			R & \text{if } q = d, \\
			0 & \text{otherwise}.
		\end{cases}
	\end{equation*}
	It follows that $ {}^{\Pi}E^2_{pq} = 0 $ for $ q \neq d $, so $ {}^{\Pi}E^2 $ has a single non-zero row given by
	\[
	{}^{\Pi}E^2_{p,d} = \overline{H}_p(\widetilde{BG}).
	\]
	In particular, the spectral sequence degenerates at the $E^2$-page, and hence $ {}^{\Pi}E^{\infty} = {}^{\Pi}E^2 $.
	Consequently,
	\[
	\overline{H}_{p+d}(\widetilde{\Pi_{hG}}) = \overline{H}_p(\widetilde{BG}),
	\]
	or equivalently $ \overline{H}_*(\widetilde{\Pi_{hG}}) = \overline{H}_{*-d}(\widetilde{BG}) $.

	\begin{table}[!htbp]
		\centering
		\caption{Upper right corner of $ {}^{\Pi}E^2 $}
		\label{table:PE2_bm_homol}
		\begin{tabular}{|r|c|c|c|c|c|c|c|c|c|}
			\hline
			$d$ &  $ 4\Q $ &       &  $ 3\Q $     &       & $2\Q $ &       & $ \Q $ &       & $ \Q $ \\ \hline
			    & $N-8$    & $N-7$ & $N-6$        & $N-5$ & $N-4$  & $N-3$ & $N-2$  & $N-1$ & $N$       \\ \hline
		\end{tabular}
	\end{table}

	The open inclusion $ \UU \hookrightarrow \Pi $ induces an open inclusion
	\[
	\UU \times_G \widetilde{EG} \hookrightarrow \Pi \times_G \widetilde{EG}.
	\]
	This yields, on the one hand, a map
	\[
	f^! \colon \overline{H}^{G}_{*}(\Pi) \to \overline{H}^{G}_{*}(\UU),
	\]
	and, on the other hand, a morphism of spectral sequences
	\[
	f^r_{pq} \colon {}^{\Pi}E^r_{p,q} \to {}^{\UU}E^r_{p,q}, \qquad r \geq 2,
	\]
	which is compatible with the differentials and therefore induces a map on the abutments.
	By \autoref{lem:open_inclusion_spse}, this induced map coincides with $f^!$.

	By definition of the homotopy quotient, the complement of $ \widetilde{\Sigma} $ in $ \widetilde{\Pi} $ is $ \widetilde{\UU} $.
	Hence the inclusion $ g \colon \Sigma \to \Pi $ induces a long exact sequence in equivariant Borel--Moore homology:
	\begin{equation*}
		\begin{tikzcd}
			\ldots \arrow[r] & \overline{H}^{G}_{n}(\Sigma) \arrow[r, "g_*"] & \overline{H}^{G}_{n}(\Pi ) \arrow[r, "f^!"] & \overline{H}^{G}_{n}( \UU ) \arrow[r] & \overline{H}^{G}_{n-1}(\Sigma) \arrow[r, "g_*"] & \ldots
		\end{tikzcd}
	\end{equation*}

	Since we work with coefficients in a field, there is a decomposition
	\[
	\overline{H}^{G}_{n}(\UU) \cong \bigoplus_p {}^{\UU}E^{\infty}_{p,n-p}.
	\]
	From the known additive structure, only the row $q = d$ is non-zero on the $E^\infty$-page, so
	\[
	\overline{H}^{G}_{n}(\UU) = {}^{\UU}E^{\infty}_{n-d,d}, \qquad
	\overline{H}^{G}_{n}(\Pi) = {}^{\Pi}E^{\infty}_{n-d,d} .
	\]

	Let $ x \in \overline{H}^{G}_{n}(\Pi) = {}^{\Pi}E^{\infty}_{n-d,d} $.
	By exactness, $ f^!(x) = 0 $ if and only if $ x = g_*(y) $ for some $ y \in \overline{H}^{G}_{n}(\Sigma) $.
	Since the spectral sequence for $\Pi$ degenerates at the $E^2$-page, we may regard $ x $ as an element of $ {}^{\Pi}E^{2}_{n-d,d} $.

	If $ f^!(x) = 0 $, then the image of the class $ x $ vanishes in $ {}^{\UU}E^{\infty}_{n-d,d} $.
	Therefore, its representative \mbox{$ f^2_{n-d,d}(x) \in {}^{\UU}E^{2}_{n-d,d} $} must be killed by a differential, and hence lies in the image of a differential.
	This proves the claim.
\end{proof}

For example, suppose that $ x = g_*(y) $ is such that $ f^6_{N-6,d}(x) \in {}^{\UU}E^{6}_{N-6,d} $ is non-zero.
Then, using the additive structure, we may conclude
\begin{equation*}
	\overline{H}^{G}_{N-6+d}( \UU ) = {}^{\UU}E^{\infty}_{N-6,d} = {}^{\UU}E^{7}_{N-6,d} = {}^{\UU}E^{6}_{N-6,d} / \im (\dif_6) = {}^{\UU}E^{6}_{N-6,d} / f^6_{N-6,d}(x) .
\end{equation*}
We find such $ x $ using the construction in the next section.

\subsection{Auxiliary spaces and integration along the fibre}
\label{sec:simplifying_to_int_along_fibre}

In this section, we construct a commutative diagram that allows us to describe certain elements $ x \in \overline{H}^{G}_{n}( \Pi ) $, as introduced at the end of the previous section, as images of an integration along the fibre map.
More precisely, they arise as the image of the ideal generated by the equivariant Euler class of a fibre bundle constructed below.
In the next section, we show that, for dimension reasons, all such elements $ x $ can be described in this way.

We let $ B $ denote the Borel subgroup of $ G = \GL_3(\C) $, and set $ P = G / B $.
We identify $ P $ with the incidence variety $ \left\{ (p,l) \in \C P^2 \times \C P^{2\vee} \mid p \in l \right\} $ of points of $ \C P^2 $ and lines through them.
For $ F \in \Pi $, let $ \Sing(F) \subset P $ be the set of pairs $ (p,l) $ such that either $ p $ is a triple point of the quartic $ \V(F) $ and $ l $ is a $4$-fold tangent to it at $ p $, or $ p $ is a tacnode of $ \V(F) $ and $ l $ is its $3$-fold tangent \cite[Section~3.1]{Dol1}.
By the description of $ \Sigma $ in the Introduction, $ F \in \Sigma $ if and only if $ \Sing(F) \neq \emptyset $.
We also denote \mbox{$ E \coloneq \left\{ (f,x) \in \Pi \times P \mid x \in \Sing(f) \right\} $}.

Consider the commutative diagram
\begin{equation*}
	\begin{tikzcd}
		E \arrow[r, two heads] \arrow[d, hook] & \Sigma \arrow[d, "g", hook] \\
		\Pi \times P \arrow[r, two heads]      & \Pi
	\end{tikzcd}
\end{equation*}
where the vertical maps are inclusions, the bottom horizontal map is a projection, and the top horizontal map is its restriction.
It is straightforward to verify that all these maps are proper.
In particular, the left vertical map is proper, as it corresponds to a morphism of vector bundles, and the total space of a sub-bundle is closed in the total space of the ambient bundle.
The maps are also equivariant.
These properties yield the induced commutative diagram on equivariant Borel--Moore homology groups.
\begin{equation*}
	\begin{tikzcd}
		\overline{H}^G_*(E) \arrow[r] \arrow[d, "\widetilde{g}"] \arrow[rd, dashed] & \overline{H}^G_*(\Sigma) \arrow[d, "g_*"] \\
		\overline{H}^G_*(\Pi  \times P) \arrow[r, "\psi"]                       & \overline{H}^G_*(\Pi )
	\end{tikzcd}
\end{equation*}
In particular, any class in $ \overline{H}^G_*(\Pi ) $ lying in the image of $ \overline{H}^G_*(E) $ also lies in the image of $ g_* $.

We now incorporate this map into a commutative diagram that enables us to determine the image of $ \overline{H}^G_*(E) $ explicitly.

Let $ \overline{E} \coloneq (\Pi  \times P) / E $ denote the quotient bundle.
Consider the following commutative diagram.
\begin{equation*}
	\begin{tikzcd}[column sep=small]
		&  & \overline{H}^G_{*}(\Pi  \times P) \arrow[dd, "Th" {inner sep=2mm}, "\sim" {rotate=90, anchor=north, inner sep=.5mm}] \arrow[rr, "\psi"]                                  &  & \overline{H}^G_*(\Pi ) \arrow[dd, "Th" {inner sep=2mm}, "\sim" {rotate=90, anchor=north, inner sep=.5mm}]                             \\
		&  &                                                                                                        &  &                                                                     \\
		\overline{H}^G_*(E) \arrow[rruu, "\widetilde{g}"] \arrow[rrdd, dashed] &  & \overline{H}^G_{*-\rk ((\Pi \times P)_{hG})}(P) \arrow[dd, "PD^{-1}" {inner sep=2mm}, "\sim" {rotate=90, anchor=north, inner sep=.5mm}] \arrow[rr]                  &  & {\overline{H}_{*-\rk (\Pi _{hG})}(\widetilde{BG})} \arrow[dd, "PD^{-1}" {inner sep=2mm}, "\sim" {rotate=90, anchor=north, inner sep=.5mm}] \\
		&  &                                                                                                        &  &                                                                     \\
		&  & {H^{\dim((P \times \widetilde{EG})/G) + \rk ((\Pi \times P)_{hG}) - *}(P_{hG})} \arrow[rr, "\int_{fibre}", dashed] &  & {H^{\dim(\widetilde{BG}) + \rk (\Pi _{hG}) - *}(BG)}
	\end{tikzcd}
\end{equation*}
Note that Poincar\'e duality holds only for sufficiently large values of $ * $.

Here, the top square is obtained by considering the diagram
\begin{equation*}
	\begin{tikzcd}
		\Pi \times P \arrow[d] \arrow[r] & \Pi \arrow[d] \\
		P \arrow[r]                      & pt
	\end{tikzcd}
\end{equation*}
where the vertical maps are $ G $-equivariant fibrations.
By taking the corresponding ``equivariant'' (in our sense) BM homology and using Thom isomorphism, we obtain the upper commutative square.
Here we used equalities $ \rk ((\Pi \times P)_{hG}) = \rk (\Pi _{hG}) = 2 \dim_{\C}(\Pi ) = d $.

Observe that $ \dim((P \times \widetilde{EG})/G) - \dim(\widetilde{BG}) = \dim(P) $.

As in \cite[Section~8.6]{Dol1}, the image of the left dashed line is given by the cup product of elements in $ H^*(P_{hG}) $ with the equivariant Euler class $e_G(\overline{E})$ of the quotient bundle.
That is, it is the ideal $ (e_G(\overline{E})) $.
By definition, the right dashed line is given by integration along the fibre.
Here we use that this map is well defined for finite-dimensional models, which are compact.
Consequently, their singular homology is naturally isomorphic to Borel--Moore homology.
See \autoref{sec:def_integration} for further details.

Recall that $ N = \dim(\widetilde{BG}) $.
Set $ * = N-p+d $.
The diagram above takes the following form:
\begin{equation*}
	\begin{tikzcd}
		& \overline{H}^G_{N-p+d}(\Pi  \times P) \arrow[dd, "Th" {inner sep=2mm}, "\sim" {rotate=90, anchor=north, inner sep=.5mm}] \arrow[r, "\psi"] & \overline{H}^G_{N-p+d}(\Pi ) \arrow[dd, "Th" {inner sep=2mm}, "\sim" {rotate=90, anchor=north, inner sep=.5mm}] \arrow[r, equal]    & {{}^{\Pi}E^2_{N-p,d}} \arrow[dd, "{f^2_{N-p,d}}" {inner sep=2mm}, "\sim" {rotate=90, anchor=north, inner sep=.5mm}] \\
		&                                                                          &                                                                                 &                                                                          &                            \\
		\overline{H}^G_{N-p+d}(E) \arrow[ruu, "\widetilde{g}"] \arrow[rdd, dashed] & \overline{H}^G_{N-p}(P) \arrow[dd, "PD^{-1}" {inner sep=2mm}, "\sim" {rotate=90, anchor=north, inner sep=.5mm}] \arrow[r]                  & {\overline{H}_{N-p}(\widetilde{BG})} \arrow[dd, "PD^{-1}" {inner sep=2mm}, "\sim" {rotate=90, anchor=north, inner sep=.5mm}] \arrow[r, equal] & {{}^{\UU}E^2_{N-p, d}} \arrow[dd, "\sim" {rotate=90, anchor=north, inner sep=.5mm}]                                          &                            \\
		&                                                                          &                                                                                 &                                                                          &                            \\
		& H^{\dim(P)+p}(P_{hG}) \arrow[r, "\int_{fibre}", dashed]          & H^{p}(BG) \arrow[r, equal]                                        & {}^{\UU}E_2^{p0}                                                           &
	\end{tikzcd}
\end{equation*}

The above results can be summarised as follows.
\begin{proposition}
	Let $ E \subset \Pi \times P $ be the subvariety defined above, and let $ \overline{E} = (\Pi \times P)/E $ denote the corresponding quotient bundle.

	Under the identifications arising from Thom isomorphism and Poincar\'e duality, the image of the map
	$
	\overline{H}^G_*(E) \longrightarrow \overline{H}^G_*(\Pi)
	$
	is contained in the image of the integration along the fibre map
	\mbox{$	H^*(P_{hG}) \longrightarrow H^*(BG)$},
	and corresponds to the image of the ideal generated by the equivariant Euler class $ e_G(\overline{E}) \in H^*_G(P) $.
\end{proposition}
We have therefore reduced the problem to the following steps.
\begin{enumerate}
	\item Calculate $ e_G = e_G(\overline E) $ explicitly.
	\item Find the ideal $ (e_G) $ it spans in $ H^*_G(P) $.
	\item Compute the image of $ (e_G) $ in $ H^*(\widetilde{BG}) $.
	\item Take the corresponding quotients by the image.
\end{enumerate}

\subsection{Explicit computation of the Euler class and the cohomology ring}

In this section, we first describe how to compute explicitly the image of integration along the fibre.
We then provide a computer algorithm that computes the image of a given polynomial under integration along the fibre map.
Next, we compute the equivariant Euler class of $ \overline E $ and the image of the ideal it generates.
Finally, we determine the cohomology ring in question.

We begin with the following general result.
See \autoref{sec:def_integration} for details.

\begin{remark}
	When the fibre bundle
	\begin{tikzcd}
		F \arrow[r] & E \arrow[r, "\pi"] & B
	\end{tikzcd}
	satisfies the assumptions of the Leray--Hirsch Theorem (LHT; see \autoref{thm:LHT}) and $ F, E, B $ are closed oriented manifolds, there exists a basis $ e_1, \ldots, e_m $ of the free $ H^*(B) $-module $ H^*(E) $, ordered by degree, such that integration along the fibre map $ \int_{fibre} = \pi_! $ is given by
	\begin{equation*}
		\pi_! \colon H^*(E) \to H^{*-\rk E}(B) \colon b_1 e_1 + \ldots + b_m e_m \mapsto b_m ,
	\end{equation*}
	where $ e_m $ is of the highest degree.
	The formula is independent of the choice of $ \{ e_i \} $ with $ i < m $.
\end{remark}

We now explain how this applies in our setting.
For $ G = \GL_3(\C) $, the Borel subgroup \mbox{$ B \sim T = \diag_3(\C) $} is homotopy equivalent to the maximal torus $ T $ -- the diagonal matrices.
Recall that $ P = G / B $, so \mbox{$ P_{hG} \coloneq ((G / B) \times EG)/G = EG/B = BB $}.
The fibre bundle in the following theorem realises the required integration along the fibre.

\begin{theorem}
	The fibre bundle
	\begin{equation*}
		\begin{tikzcd}
			G/T \arrow[r, "\iota", hook]    & BT \arrow[r, "\pi", two heads]     & BG
		\end{tikzcd}
	\end{equation*}
	satisfies the conditions of the Leray--Hirsch Theorem.
\end{theorem}

Before proving this, we introduce some notation.
Let
\begin{equation*}
	s_1 \coloneq \varepsilon_0+\varepsilon_1+\varepsilon_2,\quad s_2 \coloneq \varepsilon_0 \varepsilon_1+\varepsilon_0\varepsilon_2+\varepsilon_1\varepsilon_2,\quad s_3 \coloneq \varepsilon_0\varepsilon_1\varepsilon_2
\end{equation*}
be the elementary symmetric polynomials in the variables $ \varepsilon_0,\varepsilon_1,\varepsilon_2 $.

Then the ring
$
\Q[s_1,s_2,s_3] = \Q[\varepsilon_0,\varepsilon_1,\varepsilon_2]^{S_3}
$
is the ring of all symmetric polynomials in variables $ \varepsilon_i $ over $ \Q $, where $ S_3 $ acts by permutations of the variables $ \varepsilon_i $.

\begin{proof}
	From the following diagram and \cite[Proposition 1.3]{BGG73}, it follows that the induced map $ \iota^* $ is surjective (as $ \ker \iota^* $ is precisely the ideal $ (s_1,s_2,s_3) $).
	\begin{equation*}
		\begin{tikzcd}
			H^*(G/T) \arrow[d, equal] \arrow[r, "s"', hook, shift right] & H^*(BT) \arrow[r, "\pi_!"', shift right] \arrow[d, equal] \arrow[l, "\iota^*"', two heads, shift right] & H^{*-\dim (P)}(BG) \arrow[d, equal] \arrow[l, "\pi^*"', shift right] \\
			{\Q[\varepsilon_0,\varepsilon_1,\varepsilon_2]/(s_1,s_2,s_3)}      & {\Q[\varepsilon_0,\varepsilon_1,\varepsilon_2]}                                                                  & {\Q[\varepsilon_0,\varepsilon_1,\varepsilon_2]^{S_3}}
		\end{tikzcd}
	\end{equation*}
	Here, the vertical isomorphisms are well known.
	In particular, the middle isomorphism is obtained from the correspondence
	\begin{equation*}
		\mathfrak{X}(T) \ni ( \chi \colon T \to \C^{\times} ) \mapsto c_1(ET \times_T \C) \in H^2(BT; \Z) .
	\end{equation*}
	The left isomorphism is the one already used above: $ \iota^* $ is surjective with kernel $ (s_1, s_2, s_3) $ by \cite[Proposition 1.3]{BGG73}, so it identifies $ H^*(G/T) $ with $ \Q[\varepsilon_0,\varepsilon_1,\varepsilon_2]/(s_1,s_2,s_3) $.
	The right one identifies $ H^*(BG) $ with the ring of $ S_3 $-invariants $ \Q[\varepsilon_0,\varepsilon_1,\varepsilon_2]^{S_3} = \Q[s_1,s_2,s_3] $, as recalled in the proof of \autoref{thm:cohom_sp_of_moduli} (\cite[\textsection 1]{BGG73}).

	Therefore, the conditions of the Leray--Hirsch Theorem hold.
\end{proof}

Thus, if we construct a basis $ e_1, \ldots, e_m \in H^*(BT) $ of $ H^*(BT) $ as in the Leray--Hirsch theorem, and an algorithm to decompose each $ x \in H^*(BT) $ with respect to this basis, we will be able to compute $ \int_{fibre}(x) = \pi_!(x) $ using \autoref{prop:int_fibre_explicitly}.

To find the basis, we will use the following two lemmas.
\begin{lemma}
	\label{lem:quotient_basis_vs_ideal}
	Let $ R \coloneq \mathscr{k}[x_1, \ldots, x_n] $ be a polynomial ring over a field $ \mathscr{k} $.
	Let $ I \subset R $ be a homogeneous ideal (i.e. an ideal such that if $ \sum f_k = f \in I $, $ f_k $ homogeneous of degree $ k $, then $ f_k \in I $ for each $ k $).
	Let $ f_1 , \ldots , f_k \in R $ be homogeneous of degree $ d $, and let \mbox{$ \overline f_i \coloneq f_i \mod I $}, $ i = 1, \ldots, k $.
	Then $ \overline f_1 , \ldots , \overline f_k \in R / I $ are linearly independent over $ \mathscr{k} $ if and only if \mbox{$ f_i \notin (I, f_1, \ldots, f_{i-1}) $} for each $ i = 1, \ldots , k $.
\end{lemma}

\begin{proof}
	Suppose that $ \overline f_1 , \ldots , \overline f_k \in R / I $ are linearly dependent: $ a_1 \overline f_{1} + \ldots + a_k \overline f_{k} = 0 $.
	Then $ a_1 f_{1} + \ldots + a_k f_{k} = a_0 g \in I $.
	Let $ i \coloneq \max \{ j \mid a_j \neq 0 \} $.
	Then $ f_i = a'_0 g + a'_1 f_1 + \ldots + a'_{i-1} f_{i-1} $, where $ a'_j \coloneq - a_j / a_i $.
	Thus $ f_i \in (I, f_1, \ldots, f_{i-1}) $.

	Suppose that $ \overline f_1 , \ldots , \overline f_k \in R / I $ are linearly independent.

	Assume that $ f_i \in (I, f_1, \ldots, f_{i-1}) $ for some $ i $.
	Note that $ f_i \in (I, f_1, \ldots, f_{i-1}) $ if and only if  $ f_i = \alpha g + \beta_1 f_1 + \ldots + \beta_{i-1} f_{i-1} $ for some $ \alpha, \beta \in R $, $ g \in I $.
	In this case, since $ f_i $ is of degree $ d $, terms on the right that are not of degree $ d $ cancel out, and we have $ f_i = \tilde \alpha \tilde g + b_1 f_1 + \ldots + b_{i-1} f_{i-1} $, where $ b_j \in \mathscr{k} $ is the constant term of $ \beta_j $, and $ \deg \tilde \alpha \tilde g = d $.

	We show that $ \tilde \alpha \tilde g \in I $, arriving at contradiction, as in this case there is a linear dependence \mbox{$ \overline f_i =  b_1 \overline f_1 + \ldots + b_{i-1} \overline f_{i-1} $}.
	Note that if $ g = g_0 + g_1 + \ldots + g_m \in I $, where $ g_j $ is homogeneous of degree $ j $, then also $ g_j \in I $, as $ I $ is homogeneous.
	But then $ \tilde \alpha \tilde g = \sum \alpha_{d-j} g_j $, where $ \alpha_{d-j} \in R $ is homogeneous of degree $ d - j $.

	The contradiction shows that $ f_i \notin (I, f_1, \ldots, f_{i-1}) $ for all $ i $.
\end{proof}

For our case, we are interested in $ \mathscr{k} = \Q $, $ R = \Q[\varepsilon_0, \varepsilon_1, \varepsilon_2] $, $ I = (s_1,s_2,s_3) $.

\begin{lemma}
	\label{lem:quotient_basis_vs_module_basis}
	A family $ \{ e_i \mid i = 1, \ldots, m \} $ is a basis of the module $ \Q[\varepsilon_0,\varepsilon_1,\varepsilon_2] $ over the ring $ \Q[s_1,s_2,s_3] $ if and only if the tuple $ \{ e_i \mod (s_1, s_2, s_3) \mid i = 1, \ldots, m \} $ is a $ \Q $-basis of $ \Q[\varepsilon_0,\varepsilon_1,\varepsilon_2] / (s_1, s_2, s_3) $.
\end{lemma}

\begin{proof}
	If $\{e_i\}$ is a basis of $R$ over $\Q[s_1,s_2,s_3]$, then reducing modulo $(s_1,s_2,s_3)$ shows that their images span $R/I$ over $\Q$.
	Linear independence follows from \autoref{lem:quotient_basis_vs_ideal}.

	Conversely, suppose that the images of $\{e_i\}$ form a $\Q$-basis of $R/I$.
	Then every $f \in R$ can be written as
	\[
	f \equiv \sum_i a_i e_i \mod I, \qquad a_i \in \Q,
	\]
	so that
	\[
	f = \sum_i a_i e_i + h, \quad h \in I.
	\]
	Since $I = (s_1,s_2,s_3)$, we can write $h = \sum_j s_j h_j$.
	If $f \neq 0$, then $\deg h_j < \deg f$.
	Iterating this procedure for each $h_j$, we eventually get an expression of $f$ as an $\Q[s_1,s_2,s_3]$-linear combination of the $e_i$.

	Uniqueness follows from linear independence modulo $I$.

	Alternatively, this result follows from the Leray--Hirsch Theorem.
\end{proof}

Let $ e_1, \ldots, e_6 $ be the elements $ 1, \varepsilon_1, \varepsilon_2, \varepsilon_1^2, \varepsilon_2^2, \varepsilon_1 \varepsilon_2^2 \in \Q[\varepsilon_0,\varepsilon_1,\varepsilon_2] $.
They are homogeneous, of degrees $ 0, 1, 1, 2, 2, 3 $, and $ I = (s_1,s_2,s_3) $ is a homogeneous ideal, so $ R/I $ is graded and $ \Q $-linear independence of their images may be tested one degree at a time.
The hypothesis of \autoref{lem:quotient_basis_vs_ideal} is that the elements share a degree, which the six taken together do not.
Within a single degree the lemma does apply, and its criterion amounts there to
\begin{equation*}
	1 \notin I, \quad \varepsilon_1 \notin I, \quad \varepsilon_2 \notin (I, \varepsilon_1), \quad \varepsilon_1^2 \notin I, \quad \varepsilon_2^2 \notin (I, \varepsilon_1^2), \quad \varepsilon_1 \varepsilon_2^2 \notin I .
\end{equation*}
A direct computation verifies these (see \cite[linear\_independence.sage]{github_calc}), so the images of $ e_1, \ldots, e_6 $ are linearly independent in $ R/I $.

From the cohomology spectral sequence associated with
\begin{equation*}
	\begin{tikzcd}
		\C P^1 \arrow[r] & \PP T \C P^2 \arrow[r] & \C P^2
	\end{tikzcd}
\end{equation*}
we obtain
\begin{gather*}
	\Q[\varepsilon_0,\varepsilon_1,\varepsilon_2] / (s_1, s_2, s_3) = H^*(\GL_3(\C) / B; \Q) = H^*(\PP T \C P^2; \Q) = \Q[\overset{2}{u},\overset{2}{v}] / (u^2, v^3) = \\
	= [
	\overset{0}{\Q},
	\overset{2}{2\Q},
	\overset{4}{2\Q},
	\overset{6}{\Q}
	] .
\end{gather*}

Therefore, the images of $ e_1, \ldots, e_6 $ span the entire $ \Q $-vector space $ \Q[\varepsilon_0, \varepsilon_1, \varepsilon_2] / (s_1, s_2, s_3) $.
But then by \autoref{lem:quotient_basis_vs_module_basis} they form a basis of the $ \Q[s_1, s_2, s_3] $-module $ \Q[\varepsilon_0,\varepsilon_1,\varepsilon_2] $:
for any \mbox{$ f \in \Q[\varepsilon_0,\varepsilon_1,\varepsilon_2] $}, there is a decomposition $ f = f_1 e_1 + \ldots + f_6 e_6 $ for some symmetric polynomials \mbox{$ f_i = f_i(\varepsilon_0,\varepsilon_1,\varepsilon_2) \in \Q[s_1,s_2,s_3] $}.

Given a polynomial $ f \in \Q[\varepsilon_0,\varepsilon_1,\varepsilon_2] $, we determine its decomposition with respect to the basis $ \{ e_i \} $ using SageMath again.
More precisely, we use the algorithm described below (which is also implemented in the file \cite[decomposition.sage]{github_calc}).

Before proceeding, we require the following technical lemma from computational commutative algebra.
\begin{lemma}
	\label{lem:normal_form}
	Fix a monomial ordering on $ R $ and let $ I \subset R $ be an ideal.
	Every $ f \in R $ has a unique \emph{normal form} $ \nor(f) \in R $ modulo $ I $: a polynomial such that $ f - \nor(f) \in I $ and no monomial of $ \nor(f) $ lies in the initial ideal of $ I $.
\end{lemma}
This follows from \cite[Theorem 1.6.7]{AdamsLoust1994}.

The normal form is not, in general, a $ \Q $-linear combination of $ e_1, \ldots, e_m $: its monomials are the standard monomials of $ I $, and those form a divisor-closed set, whereas $ \{ 1, \varepsilon_1, \varepsilon_2, \varepsilon_1^2, \varepsilon_2^2, \varepsilon_1\varepsilon_2^2 \} $ contains $ \varepsilon_1\varepsilon_2^2 $ but not $ \varepsilon_1\varepsilon_2 $.
Passing from one to the other is a change of basis in $ R/I $, and is the second half of step $ 1 $ below.

\begin{algorithm}[Decomposition with respect to a Leray--Hirsch basis]
	\label{alg:decomposition_LH}
	Let $ R = \Q[\varepsilon_0,\varepsilon_1,\varepsilon_2] $ and let $ I = (s_1,s_2,s_3) \subset R $ be the ideal generated by the elementary symmetric polynomials.
	Let $ \{ e_1,\ldots,e_m \} \subset R $ be a basis of $ R $ as a free module over $ \Q[s_1,s_2,s_3] $.

	For a given polynomial $ f \in R $, the following procedure computes its decomposition
	\[
	f = \sum_{i=1}^m f_i \, e_i, \qquad f_i \in \Q[s_1,s_2,s_3].
	\]

	\begin{enumerate}
		\item Compute the normal form $ \nor(f) $ of $ f $ with respect to the ideal $ I $, and then, by linear algebra in the finite-dimensional $ \Q $-vector space $ R/I $, the coefficients $ a_1, \ldots, a_m \in \Q $ with
		\[
		\nor(f) \equiv \sum_{i=1}^m a_i e_i \pmod I ,
		\]
		which exist and are unique because the $ \overline{e_i} $ form a $ \Q $-basis of $ R/I $.
		Set $ h_0 \coloneq f - \sum_{i=1}^m a_i e_i \in I $.

		\item If $ h_0 = 0 $, then the decomposition is given by
		\[
		f = \sum_{i=1}^m a_i e_i,
		\]
		and the procedure terminates.

		\item Otherwise, express $ h_0 $ in terms of the generators of $ I $:
		\[
		h_0 = h_{01} s_1 + h_{02} s_2 + h_{03} s_3,
		\]
		where $ h_{0j} \in R $.

		\item Apply the procedure recursively to each non-zero polynomial $ h_{0j} $, obtaining decompositions
		\[
		h_{0j} = \sum_{i=1}^m f_{ij} \, e_i, \qquad f_{ij} \in \Q[s_1,s_2,s_3].
		\]

		\item Substituting these expressions into the identity for $ f $ and collecting terms yields
		\[
		f = \sum_{i=1}^m f_i \, e_i, \qquad f_i \in \Q[s_1,s_2,s_3].
		\]
	\end{enumerate}
\end{algorithm}

The procedure terminates after finitely many steps, as the degrees of the polynomials decrease at each stage.
Its correctness follows from the existence and uniqueness of normal forms modulo $I$ and from the fact that $\{e_i\}$ forms a basis of $R$ as a $\Q[s_1,s_2,s_3]$-module.

For the implementation in Python using SageMath, see \cite[decomposition.sage]{github_calc}.
For more details on computational commutative algebra we refer the reader to \cite{AdamsLoust1994}.

Let us now find the equivariant Euler class $ e_G = e_G(\overline E) = e( \overline E_{hG}) \in H^*(BT) $, so that we can use it as input for the provided SageMath algorithm.
Recall that $ \overline{E} = (\Pi \times P) / E $, where \mbox{$ E = \left\{ (f,x) \in \Pi \times P \mid x \in \Sing(f) \right\} $}.

\begin{theorem}\label{thm:stable_euler_class}
	The equivariant Euler class of $ \overline{E} $ is
	\begin{equation*}
		e_G = \prod \dif \chi_i = 4\varepsilon_0(3\varepsilon_0 +\varepsilon_1)(3\varepsilon_0 +\varepsilon_2)(2\varepsilon_0 +2\varepsilon_1)(\varepsilon_0 + 3\varepsilon_1)(2\varepsilon_0 + \varepsilon_1 + \varepsilon_2) .
	\end{equation*}
\end{theorem}

\begin{proof}
	From \cite[Lemma 4.1.6]{GorKon21} we know that $ e_G = \prod \dif  \chi_i $, where $ \chi_i $ are the characters corresponding to the equivariant vector bundle $ \overline{E} $, as above.
	Let us find it explicitly.

	Let $ p = (1:0:0), q=(0:1:0) $, and let $ l =\overline{pq} $.

	Then the fibre $ E_{(p,l)} $ of $ E $ over $ (p,l) \in P $, just as in \cite[Theorem~5.1]{Dol1}, is spanned by the following monomials:
	\begin{equation*}
		\langle y^4, y^3z, y^2z^2, yz^3, z^4 ; xy^2z, xyz^2, xz^3 ; x^2z^2 \rangle .
	\end{equation*}
	Consequently, the fibre of $ \overline{E} $ is $ \overline{E}_{(p,l)} = $
	\begin{equation*}
		= \langle x^4 ; x^3y, x^3z ; x^2y^2, x^2yz ; xy^3 \rangle .
	\end{equation*}
	We now compute the characters for the action $ (g \cdot F)(v) = F(g^{-1}v) $:
	\begin{equation*}
		T(s,t,u)(x^ny^mz^{k}) = (s^{-1}x)^n (t^{-1}y)^m (u^{-1}z)^{k}= s^{-n} t^{-m} u^{-k} x^n y^m z^{k} =\chi(s,t,u) x^ny^mz^{k}
	\end{equation*}
	yields
	\begin{equation*}
		\chi(s,t,u) = s^{-n} t^{-m} u^{-k} \implies \dif \chi = -(n\varepsilon_0 + m\varepsilon_1 + k\varepsilon_2) .
	\end{equation*}

	\begin{enumerate}
		\item $ T(s,t,u)(x^4) = s^{-4}x^4 \implies \chi_1 = s^{-4} \implies \dif \chi_1 = -4\varepsilon_0 $
		\item $ T(s,t,u)(x^3y) = (s^{-1}x)^3(t^{-1}y)=s^{-3}t^{-1} x^3y \implies \chi_2 = s^{-3}t^{-1} \implies \dif \chi_2 = -(3\varepsilon_0 +\varepsilon_1) $
		\item $ T(s,t,u)(x^3z) = (s^{-1}x)^3(u^{-1}z) = s^{-3}u^{-1} x^3z \implies \chi_3 = s^{-3}u^{-1} \implies \dif \chi_3 = -(3\varepsilon_0 +\varepsilon_2) $
		\item \sloppy $ T(s,t,u)(x^2y^2) = (s^{-1}x)^2(t^{-1}y)^2 = s^{-2}t^{-2} x^2y^2 \implies \chi_4 = s^{-2}t^{-2} \implies \mbox{$ \dif \chi_4 = -(2\varepsilon_0 +2\varepsilon_1) $} $
		\item $ T(s,t,u)(xy^3) = (s^{-1}x)(t^{-1}y)^3 = s^{-1}t^{-3} xy^3 \implies \chi_5 = s^{-1}t^{-3} \implies \dif \chi_5 = -(\varepsilon_0 + 3\varepsilon_1) $
		\item $ T(s,t,u)(x^2yz) = (s^{-1}x)^2(t^{-1}y)(u^{-1}z) = s^{-2}t^{-1}u^{-1} x^2yz \implies \chi_6 = s^{-2}t^{-1}u^{-1} \implies \dif \chi_6 = -(2\varepsilon_0 + \varepsilon_1 + \varepsilon_2) $
	\end{enumerate}

	Now, the equivariant Euler class is
	\begin{equation*}
		e_G = \prod \dif \chi_i = 4\varepsilon_0(3\varepsilon_0 +\varepsilon_1)(3\varepsilon_0 +\varepsilon_2)(2\varepsilon_0 +2\varepsilon_1)(\varepsilon_0 + 3\varepsilon_1)(2\varepsilon_0 + \varepsilon_1 + \varepsilon_2) .
	\end{equation*}
\end{proof}

Since $ H^*(BT) $ is a free \mbox{$ H^*(BG) $-module} on $ e_1, \ldots, e_6 $ and the map $ \int $ is a morphism of \mbox{$ H^*(BG) $-modules} (see \autoref{sec:def_integration}), the image $ \int (e_G) $ of the ideal $ (e_G) \subset H^*(BT) $ is the ideal of \mbox{$ H^*(BG) = \Q[s_1,s_2,s_3] $} generated by the six classes $ \int e_i e_G $.
Substituting $ e_G $ into the implementation \cite[decomposition.sage]{github_calc}, we compute them, together with the redundant class \mbox{$ \int \varepsilon_0 e_G = s_1 \int e_G - \int \varepsilon_1 e_G - \int \varepsilon_2 e_G $}, which \autoref{thm:cohom_ring_moduli} uses:
\begin{alignat*}{2}
	\int e_G                               & = -8 (72 s_1^3 - 39 s_1 s_2 + 7 s_3)                               & \quad & \in H^6(BG) ,    \\
	\int \varepsilon_0 e_G                 & = -16(36 s_1^4 - 51 s_1^2 s_2 +22s_1 s_3 + 3s_2^2)                 &       & \in H^8(BG) ,    \\
	\int \varepsilon_1 e_G                 & = -16(9 s_1^4 + 12 s_1^2 s_2 - 31 s_1 s_3)                         &       & \in H^8(BG) ,    \\
	\int \varepsilon_2 e_G                 & = 8 (18 s_1^4 - 39 s_1^2 s_2 - 25 s_1 s_3 + 6 s_2^2)               &       & \in H^8(BG) ,    \\
	\int \varepsilon_1^2 e_G               & = -8 (9 s_1^3 s_2 - 49 s_1^2 s_3 + 12 s_1 s_2^2 - 12 s_2 s_3)     &       & \in H^{10}(BG) , \\
	\int \varepsilon_2^2 e_G               & = 8 (18 s_1^3 s_2 + 54 s_1^2 s_3 - 21 s_1 s_2^2 - 11 s_2 s_3)     &       & \in H^{10}(BG) , \\
	\int \varepsilon_1 \varepsilon_2^2 e_G & = 16 s_3 (9 s_1^3 - 12 s_1 s_2 - 7 s_3)                            &       & \in H^{12}(BG) .
\end{alignat*}
The three classes of degree $ 8 $ are linearly independent.
By a dimension argument, this determines $ H^6(\UU/G) $ and $ H^8(\UU/G) $: the class of degree $ 6 $ cuts the $ 3 $-dimensional degree $ 6 $ part of $ \Q[s_1,s_2,s_3] $ down to the dimension $ 2 $ found in \autoref{thm:cohom_sp_of_moduli}, and the three of degree $ 8 $ cut its $ 4 $-dimensional degree $ 8 $ part down to dimension $ 1 $.
The classes of degrees $ 10 $, $ 10 $ and $ 12 $ contribute one further relation, in degree $ 10 $ (see \autoref{cor:cohom_ring_moduli_presentation} below).

We thus obtain the rational cohomology groups $ H^j(\UU/G) $ expressed in terms of the generating polynomials $ s_i $.

\begin{remark}
	The below relations were computed as follows.
	We found each $ \int_{fibre} \varepsilon_i e_G $ via a direct computation (see \cite[decomposition.sage]{github_calc}), then put the coefficients of each resulting polynomial into rows of a matrix, and then reduced it to row-echelon form.
\end{remark}

\begin{theorem}\label{thm:cohom_ring_moduli}
	The cohomology ring $ H^*(\UU/G) $ is non-zero only in the following degrees:
	\begin{align*}
		H^0(\UU/G) & = \langle 1 \rangle ,                                                                                                                             \\
		H^2(\UU/G) & = \langle s_1 \rangle ,                                                                                                                           \\
		H^4(\UU/G) & = \langle s_1^2, s_2 \rangle ,                                                                                                                    \\
		H^6(\UU/G) & = \langle s_1^3,\ s_1s_2,\ s_3\ |\ 72s_1^3 - 39s_1s_2 + 7s_3 \rangle ,                                                                            \\
		H^8(\UU/G) & = \langle s_1^4,\ s_1^2s_2,\ s_1s_3,\ s_2^2\ |\ \int \varepsilon_0 e_G, \int \varepsilon_1 e_G, \int \varepsilon_2 e_G \rangle =                  \\
		& = \langle s_1^4,\ s_1^2s_2,\ s_1s_3,\ s_2^2\ |\ s_1^4 - \frac{25}{369}s_2^2,\ s_1^2s_2 - \frac{17}{123}s_2^2,\ s_1s_3-\frac{3}{41}s_2^2 \rangle .
	\end{align*}
\end{theorem}

Finally, the multiplication structure is induced by the multiplication rules for polynomials.
Assembling the degrees into a single presentation requires one relation beyond those listed above.

\begin{corollary}
	\label{cor:cohom_ring_moduli_presentation}
	As a graded ring, $ H^*(\UU/G) = \Q[s_1,s_2,s_3] / I $, where $ I $ is the ideal generated by
	\begin{gather*}
		72s_1^3 - 39s_1s_2 + 7s_3 \in H^6 , \\
		s_1^4 - \frac{25}{369}s_2^2 ,\quad s_1^2s_2 - \frac{17}{123}s_2^2 ,\quad s_1s_3 - \frac{3}{41}s_2^2 \in H^8 , \\
		s_1s_2^2 \in H^{10} .
	\end{gather*}
\end{corollary}

\begin{proof}
	Write $ J \coloneq \int (e_G) $, so that $ H^*(\UU/G) = \Q[s_1,s_2,s_3] / J $, and let $ I $ be the ideal generated by the five relations listed above.
	The first four of them lie in $ J $ by the computation above.
	For the fifth, multiply the three degree $ 8 $ relations by $ s_1 $, and substitute $ 7s_3 = -72s_1^3 + 39s_1s_2 $ into $ s_2s_3 $.
	Modulo the four relations of degrees $ 6 $ and $ 8 $ this gives
	\begin{equation*}
		s_1^5 \equiv \frac{25}{369} s_1s_2^2 ,\quad s_1^3s_2 \equiv \frac{17}{123} s_1s_2^2 ,\quad s_1^2s_3 \equiv \frac{3}{41} s_1s_2^2 ,\quad s_2s_3 \equiv \frac{1191}{287} s_1s_2^2 .
	\end{equation*}
	So the degree $ 10 $ part of the quotient by those four relations is spanned by $ s_1s_2^2 $.
	Since \mbox{$ H^{10}(\UU/G) = 0 $} by \autoref{thm:cohom_sp_of_moduli}, the degree $ 10 $ part of $ J $ is all of $ \Q[s_1,s_2,s_3]_{10} $.
	In particular $ s_1s_2^2 \in J $.
	Explicitly, $ \int \varepsilon_1^2 e_G \equiv \frac{92160}{287} s_1s_2^2 $ modulo the same four relations.

	Hence $ I \subseteq J $.
	A direct computation shows that $ \Q[s_1,s_2,s_3] / I $ has dimensions $ 1,1,2,2,1 $ in degrees $ 0,2,4,6,8 $ and vanishes in all higher degrees, which by \autoref{thm:cohom_sp_of_moduli} is exactly $ \dim H^*(\UU/G) $.
	Therefore $ I = J $.
\end{proof}

This completes the description of the multiplicative structure on $ H^{*}(\UU_{hG}) = H^*(\UU/G) $.

\section{Cohomology ring of the space of equations \texorpdfstring{$\UU$}{}}
\label{sec:mult_struct_eqs}

Let $ {}^{G}E_* $ be the cohomological spectral sequence of
\begin{equation*}
	G \to \UU \to \UU / G .
\end{equation*}
We use this spectral sequence together with the multiplicative structures on $ H^*(G) $ and $ H^*(\UU / G) $ to define a multiplication on $ H^*(\UU) $.
To this end, we determine explicitly how the differentials act on the groups $ {}^{G}E^{p,q}_r $, using the decomposition $ {}^{G}E_2 = H^*(G) \otimes H^*(\UU / G)$ (for details, see \cite[Corollary A.1]{BerGorKon23}) and the multiplicative structure of the spectral sequence.

For brevity, we write $ ab \coloneq a \otimes b $.

Recall from \autoref{table:GE_2} the description of the page $ {}^{G}E^{}_{2} $.
In terms of generators, this page is given in \autoref{table:GE_2_gens}.

\begin{table}[!htbp]
	\centering
	\caption{The page ${}^{G}E^{}_{2}$ in terms of generators}
	\label{table:GE_2_gens}
	\renewcommand{\arraystretch}{1.2}
	\begin{tabular}{|r|C|C|C|C|C|C|C|C|C|}
		\hline
		& 0           & 1 & 2               & 3 & 4                                  & 5 & 6                                       & 7 & 8                   \\
		\hline
		9 & \xi_1 \xi_3 \xi_5 &   & \xi_1 \xi_3 \xi_5 s_1 &   & \makecell{\xi_1 \xi_3 \xi_5 s_1^2\\ \xi_1 \xi_3 \xi_5 s_2} &   & \makecell{\xi_1 \xi_3 \xi_5[s_1^3]\\ \xi_1 \xi_3 \xi_5[s_1s_2]} &   & \xi_1 \xi_3 \xi_5 [s_2^2] \\
		\hline
		8 & \xi_3 \xi_5     &   & \xi_3 \xi_5 s_1     &   & \makecell{\xi_3 \xi_5 s_1^2\\ \xi_3 \xi_5 s_2}         &   & \makecell{\xi_3\xi_5[s_1^3]\\ \xi_3\xi_5[s_1s_2]}           &   & \xi_3 \xi_5 [s_2^2]     \\
		\hline
		7 &             &   &                 &   &                                    &   &                                         &   &                     \\
		\hline
		6 & \xi_1 \xi_5     &   & \xi_1 \xi_5 s_1     &   & \makecell{\xi_1 \xi_5 s_1^2\\ \xi_1 \xi_5 s_2}         &   & \makecell{\xi_1\xi_5[s_1^3]\\ \xi_1\xi_5[s_1s_2]}           &   & \xi_1 \xi_5 [s_2^2]     \\
		\hline
		5 & \xi_5         &   & \xi_5 s_1         &   & \makecell{\xi_5 s_1^2\\ \xi_5 s_2}                 &   & \makecell{\xi_5[s_1^3]\\ \xi_5[s_1s_2]}                 &   & \xi_5 [s_2^2]         \\
		\hline
		4 & \xi_1 \xi_3     &   & \xi_1 \xi_3 s_1     &   & \makecell{\xi_1 \xi_3 s_1^2\\ \xi_1 \xi_3 s_2}         &   & \makecell{\xi_1\xi_3[s_1^3]\\ \xi_1\xi_3[s_1s_2]}           &   & \xi_1 \xi_3 [s_2^2]     \\
		\hline
		3 & \xi_3         &   & \xi_3 s_1         &   & \makecell{\xi_3 s_1^2\\ \xi_3 s_2}                 &   & \makecell{\xi_3[s_1^3]\\ \xi_3[s_1s_2]}                 &   & \xi_3 [s_2^2]         \\
		\hline
		2 &             &   &                 &   &                                    &   &                                         &   &                     \\
		\hline
		1 & \xi_1         &   & \xi_1 s_1         &   & \makecell{\xi_1 s_1^2\\ \xi_1 s_2}                 &   & \makecell{\xi_1 [s_1^3]\\ \xi_1[s_1s_2]}                &   & \xi_1 [s_2^2]         \\
		\hline
		0 & 1           &   & s_1             &   & \makecell{s_1^2\\ s_2}                         &   & \makecell{\left[s_1^3\right]\\ \left[s_1s_2\right]}                       &   & [s_2^2]             \\
		\hline
		& 0           & 1 & 2               & 3 & 4                                  & 5 & 6                                       & 7 & 8                   \\
		\hline
	\end{tabular}
\end{table}

More precisely, we have $ H^*(G) = \Lambda[\xi_1, \xi_3, \xi_5] $, where $ \xi_i $ has degree $ i $.
From earlier, we have an explicit presentation of the ring $ H^*(\UU / G) $ in terms of the generators $ s_i $ of $ \Q[s_1, s_2, s_3] $ and the relations among them.
The notation $ \left[ f(s_1,s_2,s_3) \right] $ in \autoref{table:GE_2_gens} denotes the class of $ f(s_1,s_2,s_3) $ modulo these relations.
It follows from \cite[Corollary A.1]{BerGorKon23} that
\begin{equation*}
	{}^{G}E_2^{pq} = H^p(\UU / G; H^q(G)) = H^p(\UU / G) \otimes H^q(G) = {}^{G}E_2^{p0} \otimes {}^{G}E_2^{0q} .
\end{equation*}

\subsection{Differentials in the spectral sequence \texorpdfstring{$ {}^{G}E_* $}{}}\label{sec:eqs_diffs}

\begin{theorem}
	The page $ {}^{G}E^{}_{\infty} $ is as in \autoref{table:GE_inf_gens}.

	\begin{table}[!htbp]
		\centering
		\caption{The page ${}^{G}E^{}_{\infty}$ in terms of generators}
		\label{table:GE_inf_gens}
		\renewcommand{\arraystretch}{1.3}
		\begin{tabular}{|r|C|C|C|C|C|C|C|C|C|}
			\hline
			& 0     & 1 & 2 & 3 & 4                            & 5 & 6                                         & 7 & 8                         \\
			\hline
			9 &       &   &   &   &                              &   & \xi_1 \xi_3 \xi_5 [51 s_1^3 - 25 s_1 s_2] &   & \xi_1 \xi_3 \xi_5 [s_2^2] \\
			\hline
			8 &       &   &   &   & \xi_3 \xi_5 (s_2 \mod s_1^2) &   &                                           &   &                           \\
			\hline
			7 &       &   &   &   &                              &   &                                           &   &                           \\
			\hline
			6 &       &   &   &   &                              &   & \xi_1 \xi_5 [51 s_1^3 - 25 s_1 s_2]       &   & \xi_1 \xi_5 [s_2^2]       \\
			\hline
			5 & \xi_5 &   &   &   &                              &   &                                           &   &                           \\
			\hline
			4 &       &   &   &   &                              &   & \xi_1 \xi_3 [51 s_1^3 - 25 s_1 s_2]       &   & \xi_1 \xi_3 [s_2^2]       \\
			\hline
			3 &       &   &   &   & \xi_3 (s_2 \mod s_1^2)       &   &                                           &   &                           \\
			\hline
			2 &       &   &   &   &                              &   &                                           &   &                           \\
			\hline
			1 &       &   &   &   &                              &   & \xi_1 [51 s_1^3 - 25 s_1 s_2]             &   & \xi_1 [s_2^2]             \\
			\hline
			0 & 1     &   &   &   &                              &   &                                           &   &                           \\
			\hline
			& 0     & 1 & 2 & 3 & 4                            & 5 & 6                                         & 7 & 8                         \\
			\hline
		\end{tabular}
	\end{table}

\end{theorem}

\begin{proof}
	Recall the Leibniz rule: $ \dif (ab) = (\dif a)b + (-1)^{\deg(a)}a\dif (b) $.

	For any element $ \xi \in H^q(G) $ of the $ \nth{0} $ column and any element $ f \in H^p(\UU / G) $ of the $ \nth{0} $ row, the differential on the product $ \xi \otimes f \in {}^{G}E^{p,q}_2 $ is $ \dif_2 (\xi f) = \dif_2 (\xi) f \pm \xi \dif_2 (f) = \dif_2 (\xi) f $, so it suffices to find the differential on the $ \nth{0} $ column.

	$ H^1(\UU) = 0 $ implies that $ \dif_2(\xi_1) = \lambda s_1 $ for some non-zero constant $ \lambda \in \Q^{\times} $.

	We have $ \dif_2(\xi_3) = 0 $, so by the Leibniz rule $ \dif_2(\xi_1\xi_3) = \dif_2(\xi_1)\xi_3 \pm \xi_1\dif_2(\xi_3) = \lambda s_1\xi_3 $.

	Since $ H^5(\UU) = \Q $, we have $ \dif_2(\xi_5) = 0 $, as it generates the only group on the diagonal $ p+q = 5 $ that survives.

	Hence $ \dif_2 (\xi_1 \xi_5) = \lambda s_1 \xi_5 $ and $ \dif_2(\xi_3\xi_5) = 0 $, and therefore $ \dif_2(\xi_1\xi_3\xi_5) = \lambda s_1\xi_3\xi_5 $.

	Thus $ {}^{G}E^{}_{3} $ is as in \autoref{table:GE_3_gens}, where the class $ k \coloneq [51 s_1^3 - 25 s_1 s_2] $ is chosen so that $ \dif_2^{6,1}(\xi_1 k) = 0 $, i.e. the elements in the $ \nth{6} $ column generate the kernel of $ \dif_2 $.

	\begin{table}[!htbp]
		\centering
		\caption{The page ${}^{G}E^{}_{3}$ in terms of generators}
		\label{table:GE_3_gens}
		\renewcommand{\arraystretch}{1.3}
		\begin{tabular}{|r|C|C|C|C|C|C|C|C|C|}
			\hline
			& 0           & 1 & 2 & 3 & 4                            & 5 & 6                                         & 7 & 8                         \\
			\hline
			9 &             &   &   &   &                              &   & \xi_1 \xi_3 \xi_5 [51 s_1^3 - 25 s_1 s_2] &   & \xi_1 \xi_3 \xi_5 [s_2^2] \\
			\hline
			8 & \xi_3 \xi_5 &   &   &   & \xi_3 \xi_5 (s_2 \mod s_1^2) &   &                                           &   &                           \\
			\hline
			7 &             &   &   &   &                              &   &                                           &   &                           \\
			\hline
			6 &             &   &   &   &                              &   & \xi_1 \xi_5 [51 s_1^3 - 25 s_1 s_2]       &   & \xi_1 \xi_5 [s_2^2]       \\
			\hline
			5 & \xi_5       &   &   &   & \xi_5 (s_2 \mod s_1^2)       &   &                                           &   &                           \\
			\hline
			4 &             &   &   &   &                              &   & \xi_1 \xi_3 [51 s_1^3 - 25 s_1 s_2]       &   & \xi_1 \xi_3 [s_2^2]       \\
			\hline
			3 & \xi_3       &   &   &   & \xi_3 (s_2 \mod s_1^2)       &   &                                           &   &                           \\
			\hline
			2 &             &   &   &   &                              &   &                                           &   &                           \\
			\hline
			1 &             &   &   &   &                              &   & \xi_1 [51 s_1^3 - 25 s_1 s_2]             &   & \xi_1 [s_2^2]             \\
			\hline
			0 & 1           &   &   &   & s_2 \mod s_1^2               &   &                                           &   &                           \\
			\hline
			& 0           & 1 & 2 & 3 & 4                            & 5 & 6                                         & 7 & 8                         \\
			\hline
		\end{tabular}
	\end{table}

	Now, as $ H^4(\UU) = 0 $, we have $ \dif_4(\xi_3) = \mu [s_2] $ for some non-zero constant $ \mu \in \Q^{\times} $.
	We have $ \dif_4(\xi_5) = 0 $.
	Consequently, by the Leibniz rule, we also have $ \dif_4(\xi_3 \xi_5) = \mu [s_2] \xi_5 $.

	Since no further non-zero differentials can occur, we obtain $ {}^{G}E^{}_{5} = {}^{G}E^{}_{\infty} $, as stated.
\end{proof}

\subsection{Multiplicative structure of \texorpdfstring{$ H^{*}(\UU;\Q) $}{}}

We now use the multiplicative structure of $ {}^{G} E_{*} $ to determine the multiplication in $ H^*( \UU ) $.
For brevity, we denote $ H^k(\UU; \Q) $ by $ H^k $.
We also use the notation introduced in \autoref{table:cohom_U_gens}.

\begin{table}[!htbp]
	\centering
	\caption{Definitions of some useful elements in $ H^*(\UU) $}
	\label{table:cohom_U_gens}
	\renewcommand{\arraystretch}{1.3}
	\begin{tabular}{|c|c|}
		\hline
		Symbol & Definition                                                                                     \\
		\hline
		$x$    & $ \xi_5 \in {}^{G}E^{0,5}_{\infty} = H^{5} $                                                  \\
		\hline
		$y$    & $ \xi_1 [51 s_1^3 - 25 s_1 s_2] \in {}^{G}E^{6,1}_{\infty} \subset H^{7} $                     \\
		\hline
		$y'$   & element of $ H^{7} $ corresponding to $ \xi_3 (s_2 \mod s_1^2) \in {}^{G}E^{4,3}_{\infty} $    \\
		\hline
		$z$    & $ \xi_1 [s_2^2] \in {}^{G}E^{8,1}_{\infty} = H^{9} $                                           \\
		\hline
		$t$    & $ \xi_1 \xi_3 [51 s_1^3 - 25 s_1 s_2] \in {}^{G}E^{6,4}_{\infty} = H^{10} $                    \\
		\hline
		$e$    & an element of $ H^{12} $ with $ xe \neq 0 $ (\autoref{rem:e_normalisation})                  \\
		\hline
	\end{tabular}
\end{table}

\begin{remark}
	\label{rem:e_normalisation}
	Unlike the other five, the generator $ e $ is not pinned down by \autoref{thm:eqs_mult}.
	That theorem holds for any $ e \in H^{12} $ linearly independent with $ xy $ and $ xy' $, equivalently for any $ e $ with $ xe \neq 0 $, and $ \xi_1 \xi_3 [s_2^2] \in {}^{G}E^{8,4}_{\infty} $ is one such choice.
	Another natural candidate is a class pulled back from the space of curves $ \UU/\C^{\times} $ along the quotient map $ \UU \to \UU/\C^{\times} $.
	We do not fix a choice.
\end{remark}

\begin{theorem}\label{thm:eqs_mult}
	The ring $ H^*(\UU ) $ is generated by the elements listed in \autoref{table:cohom_U_gens}, with multiplication given by \autoref{table:cohom_U_mult_table}.

	\begin{table}[!htbp]
		\centering
		\caption{Multiplication table of $ H^*(\UU) $}
		\label{table:cohom_U_mult_table}
		\renewcommand{\arraystretch}{1.4}
		\setlength{\tabcolsep}{8pt}
		\begin{tabular}{CC|CCCCCCCCCCC}
			&     & 5    & 7  & 7   & 9  & 10 & 12 & 12 & 12  & 14 & 15 & 17 \\
			&     & x    & y  & y'  & z  & t  & e  & xy & xy' & xz & xt & xe \\
			\hline
			5  & x   &      & xy & xy' & xz & xt & xe &    &     &    &    &    \\
			7  & y   & -xy  &    &     &    &    &    &    &     &    &    &    \\
			7  & y'  & -xy' &    &     &    &    &    &    &     &    &    &    \\
			9  & z   & -xz  &    &     &    &    &    &    &     &    &    &    \\
			10 & t   & xt   &    &     &    &    &    &    &     &    &    &    \\
			12 & e   & xe   &    &     &    &    &    &    &     &    &    &    \\
			12 & xy  &      &    &     &    &    &    &    &     &    &    &    \\
			12 & xy' &      &    &     &    &    &    &    &     &    &    &    \\
			14 & xz  &      &    &     &    &    &    &    &     &    &    &    \\
			15 & xt  &      &    &     &    &    &    &    &     &    &    &    \\
			17 & xe  &      &    &     &    &    &    &    &     &    &    &
		\end{tabular}
	\end{table}

\end{theorem}

\begin{proof}

	For the product $ H^i \smile H^{j} \to H^{i+j} $, all three groups are non-zero only for the following pairs $ \{i, j\} $:
	\begin{equation*}
		\{5,5\}, \{5,7\}, \{5,9\}, \{5,10\}, \{5,12\}, \{7,7\}, \{7,10\} .
	\end{equation*}

	{\boldmath $ H^5 \smile H^{5\ } \to H^{10} $:} Observe that $ H^5 = \langle x \rangle  $.
	Since the square of an odd-degree element is zero by graded anti-commutativity, this product vanishes.

	{\boldmath $ H^5 \smile H^{7\ } \to H^{12} $:} Let $  H^5 = \langle x \rangle  $ and $ H^7 = \langle y, y' \rangle  $.
	Recall that there is a filtration \mbox{$ 0 \subset F^n_n \subset \ldots \subset F^n_1 \subset F^n_0 = H^n $} such that $ {}^{G}E^{p,n-p}_{\infty} = F^n_p / F^n_{p+1} $.
	Without loss of generality, $ y \in F^7_6 \subsetneq F^7_4 = H^7 $ and $ y' \in F^7_4 \setminus F^7_6 $.

	The product $ H^5 \times H^{7\ } \to H^{12} $ restricts to the products
	\begin{align*}
		F^5_0 \times F^7_4 & \to F^{12}_{4} , \\
		F^5_0 \times F^7_6 & \to F^{12}_{6} .
	\end{align*}
	The first restriction induces the product on quotients
	\begin{equation*}
		F^5_0/0 \times F^7_4 /F^7_5 \to F^{12}_{4} / F^{12}_{5} ,
	\end{equation*}
	where $ F^5_0 / 0 = {}^{G}E^{0,5}_{\infty} = \langle \xi_5 \rangle $, $ F^7_4 /F^7_5 = {}^{G}E^{4,3}_{\infty} = \langle \xi_3 [s_2] \rangle  $ and $ F^{12}_4 /F^{12}_5 = {}^{G}E^{4,8}_{\infty} = \langle \xi_5\xi_3 [s_2] \rangle  $.

	Now, if the product $ F^5_0 \times F^7_4 \to F^{12}_{4} $ was zero, then the induced product on quotients would also be zero.
	But we know that the latter product is non-trivial.
	Therefore, $ xy' \neq 0 $.

	A similar argument applied to the second restriction shows that $ xy \neq 0 $.

	The elements $ xy $ and $ xy' $ are linearly independent, and so are $ e, xy, xy' $ (\autoref{rem:e_normalisation}).
	Hence $ H^{12} = \langle e, xy, xy' \rangle  $.

	{\boldmath $ H^7 \smile H^{7\ } \to H^{14} $:} As in the case of $ \{5,5\} $, $ y^2 = 0 = y'^2 $.

	Since $ F^{14}_{14} / 0 = {}^{G}E^{14,0}_{\infty} = 0 $, it follows that $ F^{14}_{14} = 0 $.

	Continuing inductively, $ F^{14}_{i} / F^{14}_{i+1} = {}^{G}E^{i,14-i}_{\infty} = 0 $ implies that $ F^{14}_{i} = 0 $ for $ i = 10, \ldots, 14 $.

	As the product $ H^7 \times H^{7} \to H^{14} $ restricts to the product $ F_6^7 \times F_4^7 \to F_{10}^{14} $ and $ F_{10}^{14} = 0 $, we conclude that the product $ yy' = 0 $ is zero.

	{\boldmath $ H^5 \smile H^{9\ } \to H^{14} $:} As in the case $ \{5,7\} $.

	{\boldmath $ H^5 \smile H^{10} \to H^{15} $:} As in the case $ \{5,7\} $.

	{\boldmath $ H^5 \smile H^{12} \to H^{17} $:} As in the cases $ \{5,7\} $ and $ \{7,7\} $.

	{\boldmath $ H^7 \smile H^{10} \to H^{17} $:} As in the case $ \{7,7\} $.
\end{proof}

\section*{Acknowledgements}

MD was partially funded by the Deutsche Forschungsgemeinschaft (DFG, German Research Foundation) under Germany's Excellence Strategy -- The Berlin Mathematics Research Center MATH+ (EXC-2046/2, project ID: 390685689).

\section*{Tool use disclosure}

AI tools (specifically, Claude models Opus and Fable) were used in the preparation of this paper for proofreading and drafting some parts of the document.
We subsequently re-checked and edited these parts, and accept full responsibility for the correctness of all results presented in this paper.

Where indicated, the computer algebra computations use SageMath \cite{sage}, with the code in \cite{github_calc}.

\appendix

\section{Homotopy (Borel) quotient}

Let $ G $ be a group acting on a space $ Y $.
Consider the \emph{homotopy quotient} \mbox{$ Y_{hG} \coloneq (Y \times EG) / G $}.
Here, $ EG $ is the \emph{universal $ G $-bundle}.
The group $ G $ acts on the product via \mbox{$ g . (y,e) \coloneq (gy, eg^{-1}) $}, where the action of $ G $ on $ EG $ is free and proper, as in the definition of the \emph{classifying space $ BG = EG / G $}.

Consider the fibration $ Y_{hG} \to EG / G$.
The fibre over the class of an element $ b \in EG $ is isomorphic to $ Y / G_b = Y $, since the action of $ G $ on $ EG $ is free and hence the stabiliser $ G_b $ is trivial.

Thus we have a fibre bundle
\begin{equation*}
	Y \to Y_{hG} \to BG
\end{equation*}

If the group $ G $ is path-connected, then the associated spectral sequence $ E^{pq}_{*} \implies H^{p+q}(Y_{hG}) $ converges to the cohomology of $ Y_{hG} $, and its second page is given by $ E^{pq}_{2} = H^{p}(BG; H^{q}(Y; R)) $.

For more details, see \cite[Section A.3]{BerGorKon23}.

\section{Integration along the fibre}
\label{sec:def_integration}

Let \begin{tikzcd}
	F \arrow[r] & E \arrow[r, "\pi"] & B
\end{tikzcd} be a fibre bundle, where $ F $, $ E $ and $ B $ are oriented closed manifolds.
All (co)homology groups below are taken with coefficients in $ \Q $, although the results hold over any field.

\begin{definition}
	The \emph{integration along the fibre} map is defined, via commutativity, as the bottom horizontal map in the diagram below.
	\begin{equation*}
		\begin{tikzcd}
			H_{*-\dim(E)}(E) \arrow[rr, "\pi_*"]    &  & H_{*-\dim(E)}(B) \arrow[dd, "PD_B^{-1}"] \\
			&  &                                                           \\
			H^{*}(E) \arrow[uu, "PD_E"] \arrow[rr, "\pi_!", dashed] &  & H^{*-\rk E}(B)
		\end{tikzcd}
	\end{equation*}
\end{definition}

Recall the Leray--Hirsch Theorem (LHT).
\begin{theorem}[Leray--Hirsch]
	\label{thm:LHT}
	Let \begin{tikzcd}
		F \arrow[r] & E \arrow[r, "\pi"] & B
	\end{tikzcd} be a fibre bundle such that for each degree $ j $, the singular rational cohomology vector space $ H^j(F) $ is finite-dimensional, and the inclusion map $ \iota \colon F \to E $ induces a surjection $ \iota^* \colon H^*(E) \to H^*(F) $.
	Let $ s \colon H^*(F) \to H^*(E) $ be a section of this surjection, so that $ \iota^* \circ s = Id $.
	Then the following linear map is an isomorphism (called the \emph{Leray--Hirsch isomorphism}):
	\begin{equation*}
		H^*(F) \otimes H^*(B) \to H^*(E) \colon \alpha \otimes \beta \mapsto s (\alpha)\smallsmile \pi^*(\beta) .
	\end{equation*}
\end{theorem}

In terms of coordinates, it can be restated as follows.
Suppose that for each degree $ j $ there exist classes \mbox{$ e_{1,j}, \ldots, e_{m_j,j} $} in $ H^j(E) $ such that the restrictions $ \{ e_{i,j} |_{F_x} \} $ of $ e_{i,j} $ to each fibre $ F_x $ form the basis of $ H^j(F_x) $, $ x \in B $.
Then the map below is an isomorphism of $ H^*(B) $-modules.
\begin{equation*}
	\begin{array}{ccc}
		H^*(F)\otimes H^*(B)                               & \longrightarrow & H^*(E)                                       \\
		\sum_{i,j,k} a_{i,j,k} \iota^*(e_{i,j})\otimes b_k & \longmapsto     & \sum_{i,j,k}a_{i,j,k}e_{i,j}\wedge\pi^*(b_k)
	\end{array} ,
\end{equation*}
where $\{b_k\}$ is a basis for $H^*(B)$ and thus induces a basis $\{\iota^*(e_{i,j})\otimes b_k\}$ for $H^*(F)\otimes H^*(B)$.

\begin{proposition}
	The integration along the fibre map $ \pi_! $ is a morphism of $ H^*(B) $-modules.
	That is, for all $ \alpha \in H^*(E) $ and $ \beta \in H^*(B) $ one has
	\[
	\pi_!(\alpha \smallsmile \pi^*(\beta)) = \pi_!(\alpha) \smallsmile \beta .
	\]
\end{proposition}

\begin{proof}
	This is essentially the projection formula, and it uses only the duality isomorphisms $ PD_E $ and $ PD_B $ that already enter the definition of $ \pi_! $.

	Consider the commutative diagram defining $ \pi_! $, where it is expressed via Poincaré duality and the pushforward in homology.
	The pushforward map in homology satisfies the projection formula
	\[
	\pi_*(x \frown \pi^*(\beta)) = \pi_*(x) \frown \beta
	\]
	for all $ x \in H_*(E) $ and $ \beta \in H^*(B) $.

	Recall that $ PD(\alpha \smile \beta) = PD(\alpha) \frown \beta $ for compatible degrees.
	On the one hand,
	\[
	\pi_*(PD(\alpha) \frown \pi^*(\beta)) = \pi_*(PD(\alpha)) \frown \beta = PD(\pi_{!}(\alpha)) \frown \beta = PD(\pi_{!}(\alpha) \smile \beta)
	\]
	On the other hand,
	\[
	\pi_*(PD(\alpha) \frown \pi^*(\beta)) = \pi_*(PD(\alpha \smile \pi^*(\beta))) = PD(\pi_{!}(\alpha \smile \pi^*(\beta)))
	\]

	Comparing the two expressions, we obtain
	\[
	\pi_!(\alpha \smallsmile \pi^*(\beta)) = \pi_!(\alpha) \smallsmile \beta .
	\]
	This shows that $ \pi_! $ is a morphism of $ H^*(B) $-modules.
\end{proof}

The following result is convenient to have for explicit computations:

\begin{proposition}\label{prop:int_fibre_explicitly}
	When the Leray--Hirsch theorem holds and $ F, E, B $ are closed oriented manifolds, there is a basis $ e_1, \ldots, e_m $ of the free $ H^*(B) $-module $ H^*(E) $, ordered by degree, such that integration along the fibre map $ \pi_! $ is of the form
	\begin{equation*}
		\pi_! \colon H^*(E) \to H^{*-\rk E}(B) \colon b_1 e_1 + \ldots + b_m e_m \mapsto b_m ,
	\end{equation*}
	where $ e_m $ is of the highest degree.
	The formula is independent of the choice of $ \{ e_i \} $ with $ i < m $.
\end{proposition}

\begin{proof}
	Let the basis $ e_1, \ldots, e_m $ be as in the Leray--Hirsch theorem.

	First, if $ i < m $, then for degree reasons $ \pi_!(e_i) = 0 $,  since in this case $ (\deg(e_i) - \rk E) < 0 $, and hence the target group is zero.

	Next, consider $ e_m $.
	If $ \pi_!(e_m) $ is zero, then $ \pi_! \equiv 0 $.
	But by commutativity, this would imply $ \pi_* \equiv 0 $, which cannot be true for an induced morphism, as it is non-zero at least in dimension $ 0 $.

	Finally, as we are working over a field, the element $ e_m $ can be normalised so that $ \pi_!(e_m) = 1 $.
\end{proof}

\providecommand{\CH}{\operatorname{CH}}
\providecommand{\cl}{\operatorname{cl}}

\section{The Chow ring of the stack of GIT-stable quartics}
\label{sec:chow}

The description of $ H^*(\UU/G) $ as a quotient of $ H^*(BG) $ has an exact counterpart in intersection theory.
As in the rest of the paper, $ G = \GL_3(\C) $.
The rational Chow ring of the stack of GIT-stable plane quartics has the same generators and the same relations as $ H^*(\UU/G) $, and the cycle class map identifies the two rings (\autoref{cor:chow}).

For a smooth $ G $-variety $ X $ let $ \CH^*_G(X) $ be its $ G $-equivariant Chow ring with rational coefficients, graded by codimension.
It comes with the cycle class map \cite[Section 2.8]{EdiGra98}
\begin{equation*}
	\cl \colon \CH^i_G(X) \longrightarrow H^{2i}(X_{hG}) .
\end{equation*}
The ring $ \CH^*_G(\UU) $ is the rational Chow ring of the quotient stack $ [\UU/G] $ of GIT-stable plane quartics.

\begin{corollary}
	\label{cor:chow}
	The cycle class map
	\begin{equation*}
		\cl \colon \CH^*_G(\UU) \longrightarrow H^{2*}(\UU_{hG}) = H^{2*}(\UU/G)
	\end{equation*}
	is an isomorphism of rings.
	In particular $ \CH^*_G(\UU) = \Q[s_1,s_2,s_3]/I $, where $ s_i $ has codimension $ i $ and $ I $ is the ideal of \autoref{cor:cohom_ring_moduli_presentation}.
\end{corollary}

Let us recall the notation.
The group $ G $ contains the Borel subgroup $ B $ and the maximal torus $ T \subset B $ of diagonal matrices.
The quotient $ P = G/B $ is the variety of pairs $ (p,l) $ of a point $ p $ of $ \C P^2 $ and a line $ l $ through it.
Over $ P $ there are the vector bundles
\begin{equation*}
	E = \left\{ (f,x) \in \Pi \times P \mid x \in \Sing(f) \right\} \subset \Pi \times P
	\qquad \text{and} \qquad
	\overline{E} = (\Pi \times P) / E .
\end{equation*}
The projection of $ E $ to $ \Pi $ has image $ \Sigma $.
The equivariant Euler class of the quotient bundle is
\begin{equation*}
	e_G = e_G(\overline{E}) \in H^*(P_{hG}) = H^*(BT) = \Q[\varepsilon_0,\varepsilon_1,\varepsilon_2] ,
\end{equation*}
and it is computed in \autoref{thm:stable_euler_class}.
The image of the ideal $ (e_G) $ under the integration along the fibre map is the ideal
\begin{equation*}
	J \coloneq \int (e_G) \subset H^*(BG) = \Q[s_1,s_2,s_3] .
\end{equation*}
By \autoref{thm:cohom_ring_moduli} and the proof of \autoref{cor:cohom_ring_moduli_presentation},
\begin{equation}
	\label{eq:chow_cohomology_side}
	H^*(\UU/G) = H^*(BG) / J \qquad \text{and} \qquad J = I .
\end{equation}

On the side of Chow groups we need one more piece of notation.
For a $ G $-variety $ X $ that need not be smooth, $ \CH^G_*(X) $ denotes its equivariant Chow groups with rational coefficients, graded by dimension.

\begin{proof}[Proof (\autoref{cor:chow})]
	Let $ \pi \colon P \to pt $ be the projection, and let
	\begin{equation*}
		c \coloneq c_6^G(\overline E) \in \CH^6_G(P)
	\end{equation*}
	be the top equivariant Chern class of the quotient bundle $ \overline E $, so that $ \cl(c) = e_G $.
	The proof has three parts.
	We describe $ \CH^*_G(\UU) $ as a quotient of $ \CH^*_G(pt) $ by an ideal.
	Then we show that the cycle class map carries this ideal onto $ J $.
	Finally we compare the two quotients.

	First, let us show that
	\begin{equation}
		\label{eq:chow_quotient}
		\CH^*_G(\UU) = \CH^*_G(pt) / \pi_* \bigl( c \cdot \CH^*_G(P) \bigr) .
	\end{equation}
	We identify $ \CH^*_G(\Pi) $ with $ \CH^*_G(pt) $ by homotopy invariance.
	The localisation sequence
	\begin{equation*}
		\CH^G_*(\Sigma) \longrightarrow \CH^G_*(\Pi) \longrightarrow \CH^G_*(\UU) \longrightarrow 0
	\end{equation*}
	is exact \cite[Proposition 5]{EdiGra98}.
	So the restriction $ \CH^*_G(\Pi) \to \CH^*_G(\UU) $ is onto, and its kernel is the image of $ \CH^G_*(\Sigma) $.

	Let us describe this kernel.
	The projection $ E \to \Sigma $ is proper, surjective and equivariant.
	So, with rational coefficients, the push-forward $ \CH^G_*(E) \to \CH^G_*(\Sigma) $ is onto.
	The kernel is therefore the image of the push-forward along the composite
	\begin{equation*}
		E \overset{j}{\longrightarrow} \Pi \times P \longrightarrow \Pi
	\end{equation*}
	of the inclusion $ j $ with the projection.
	Both $ \CH^G_*(E) $ and $ \CH^G_*(\Pi \times P) $ are identified with $ \CH^*_G(P) $ by pull-back along the bundle projections to $ P $.
	Since $ E $ is a sub-bundle of $ \Pi \times P $ with quotient $ \overline E $, under these identifications $ j_* $ is multiplication by $ c $ \cite[Tag 0FAA]{Stacks}.
	The push-forward along $ \Pi \times P \to \Pi $ becomes $ \pi_* $.
	Indeed, proper push-forward commutes with flat pull-back in the Cartesian square formed by $ \Pi \times P $, $ \Pi $, $ P $ and $ pt $ \cite[Tag 02RG]{Stacks}.
	So the kernel is $ \pi_* \bigl( c \cdot \CH^*_G(P) \bigr) $, which gives \eqref{eq:chow_quotient}.

	Next, let us show that the cycle class map carries the ideal in \eqref{eq:chow_quotient} onto $ J $.
	The map $ \cl $ is a ring homomorphism that commutes with restriction to open subsets and with proper push-forward, and under it $ \pi_* $ corresponds to integration along the fibre.
	It is an isomorphism for $ pt $ and for $ P = G/B $:
	\begin{equation*}
		\begin{aligned}
			\cl \colon \CH^*_G(pt) & \overset{\sim}{\longrightarrow} H^*(BG) = \Q[s_1,s_2,s_3] , \\
			\cl \colon \CH^*_G(P) = \CH^*_T(pt) & \overset{\sim}{\longrightarrow} H^*(BT) = \Q[\varepsilon_0,\varepsilon_1,\varepsilon_2] .
		\end{aligned}
	\end{equation*}
	Indeed, $ \CH^*_G(pt) $ and $ \CH^*_T(pt) $ are polynomial rings on Chern classes \cite[Section 3.2]{EdiGra98}.
	The rings $ H^*(BG) $ and $ H^*(BT) $ are polynomial rings on the topological Chern classes of the same bundles, and $ \cl $ sends Chern classes to Chern classes.
	Finally, $ \CH^*_G(G/B) = \CH^*_B(pt) = \CH^*_T(pt) $, the second equality by homotopy invariance, because $ B/T $ is an affine space.
	Since $ \cl(c) = e_G $, it follows that
	\begin{equation*}
		\cl \Bigl( \pi_* \bigl( c \cdot \CH^*_G(P) \bigr) \Bigr) = \int \bigl( e_G \cdot H^*(BT) \bigr) = J .
	\end{equation*}

	It remains to compare the two quotients.
	Consider the square
	\begin{equation*}
		\begin{tikzcd}
			\CH^*_G(pt) = \CH^*_G(\Pi) \arrow[r, two heads] \arrow[d, "\cl" {inner sep=2mm}, "\sim" {rotate=90, anchor=north, inner sep=.5mm}] & \CH^*_G(\UU) \arrow[d, "\cl" {inner sep=2mm}] \\
			H^*(BG) = H^*(\Pi_{hG}) \arrow[r, two heads]                             & H^*(\UU_{hG}) = H^*(\UU/G)
		\end{tikzcd}
	\end{equation*}
	whose horizontal maps are the restrictions from $ \Pi $ to $ \UU $.
	It commutes, because $ \cl $ commutes with restriction to open subsets.
	The kernel of the upper map is the ideal in \eqref{eq:chow_quotient}, and the kernel of the lower map is $ J $ by \eqref{eq:chow_cohomology_side}.
	The left map is an isomorphism, and it carries the first kernel onto the second.
	So the right map is an isomorphism of rings.
	The presentation $ \CH^*_G(\UU) = \Q[s_1,s_2,s_3]/I $ follows, because $ J = I $ by \eqref{eq:chow_cohomology_side}.
\end{proof}

The action of $ G $ on $ \UU $ is proper with geometric quotient $ \UU/G $.
So $ \CH^*_G(\UU) $ is also the rational Chow ring of the variety $ \UU/G $ \cite[Theorems 3 and 4]{EdiGra98}.

The rational Chow ring of the stack of GIT-stable plane quartics was computed independently by Ascher and Lee.
They present it as the quotient of $ \Q[H, c_2, c_3] $ by four relations, of degrees $ 3 $, $ 4 $, $ 4 $ and $ 5 $ \cite[Theorem 9.15]{AscLee26}.
The two presentations agree.
Indeed, the substitution
\begin{equation*}
	s_1 = 3u, \qquad s_2 = c_2 + 3u^2, \qquad s_3 = c_3 + c_2 u + u^3, \qquad H = -4u
\end{equation*}
identifies $ \Q[s_1,s_2,s_3] $ with $ \Q[H, c_2, c_3] $, and a computation with Gr\"obner bases shows that it carries the ideal $ I $ onto the ideal generated by their relations.
Under this substitution, their relation of degree $ 3 $ equals $ -8(72s_1^3 - 39s_1s_2 + 7s_3) $, which is the image of the Euler class $ e_G $ under integration along the fibre.

\section{Koszul homology of the cohomology ring of the moduli space}
\label{sec:koszul}

This appendix records a few observations that are not used elsewhere in the paper.
The Koszul homology of the ring $ H^*(\UU/G) $ over $ H^*(BG) $ coincides with the cohomology of the space of equations $ \UU $.
The dimensions are those of \autoref{thm:cohom_of_equations}, with the weights, and the products are those of \autoref{thm:eqs_mult}.
This is \autoref{prop:koszul_comparison}.

We write
\begin{equation*}
	S \coloneq H^*(BG) = \Q[s_1,s_2,s_3]
	\qquad \text{and} \qquad
	R \coloneq H^*(\UU/G) = S/I ,
\end{equation*}
where $ I $ is the ideal of \autoref{cor:cohom_ring_moduli_presentation}.
By \autoref{thm:cohom_ring_moduli}, the ring $ R $ has dimension $ 7 $ over $ \Q $, with the basis
\begin{equation*}
	1, \qquad s_1, \qquad s_1^2, \ s_2, \qquad s_1^3, \ s_1s_2, \qquad s_2^2
\end{equation*}
in degrees $ 0 $, $ 2 $, $ 4 $, $ 6 $ and $ 8 $.

Let $ \sigma_1 $, $ \sigma_2 $, $ \sigma_3 $ be three symbols, one for each generator of $ S $.
The Koszul complex of $ s_1, s_2, s_3 $ on $ R $ is the exterior algebra over $ R $ on these symbols,
\begin{equation*}
	K \coloneq R \otimes \Lambda(\sigma_1, \sigma_2, \sigma_3) .
\end{equation*}
Its differential $ d $ is $ R $-linear, sends $ \sigma_k $ to the image of $ s_k $ in $ R $, and satisfies the Leibniz rule.
For example,
\begin{equation*}
	d(\sigma_1 \sigma_2) = s_1 \sigma_2 - s_2 \sigma_1 .
\end{equation*}
The complex $ K $ has two gradings.
Let $ r \in R $ be homogeneous of degree $ 2a $, and let $ k_1 < \dots < k_i $.
The element $ r \sigma_{k_1} \cdots \sigma_{k_i} $ of $ K $ has \emph{homological degree}  $i$ and \emph{internal degree} $2a + 2(k_1 + \dots + k_i)$.
So the homological degree counts the symbols.
The internal degree is the degree of the coefficient in $ R $ plus $ 2k $ for each symbol $ \sigma_k $, where $ 2k $ is the degree of $ s_k $.

The differential lowers the homological degree by one and preserves the internal degree.
So $ K $ is a differential graded algebra, its homology is an algebra, and both gradings pass to the homology.

\begin{lemma}
	\label{lem:koszul_tor}
	The homology of $ K $ is $ \operatorname{Tor}^S(\Q, R) $.
\end{lemma}

\begin{proof}
	The elements $ s_1, s_2, s_3 $ form a regular sequence in $ S $.
	So the Koszul complex \mbox{$ S \otimes \Lambda(\sigma_1, \sigma_2, \sigma_3) $} is a free resolution of the $ S $-module $ S/(s_1,s_2,s_3) = \Q $.
	The complex $ K $ is its tensor product with $ R $ over $ S $.
\end{proof}

We write $ \operatorname{Tor}^S_i(\Q, R)_{2j} $ for the part of the homology of homological degree $ i $ and internal degree $ 2j $.

\begin{proposition}
	\label{prop:koszul_dimensions}
	The total dimension of $ \operatorname{Tor}^S(\Q, R) $ is $ 12 $, and the internal degrees of a basis are
	\begin{equation*}
		\begin{aligned}
			\operatorname{Tor}^S_0(\Q, R) & \colon \quad 0 , \\
			\operatorname{Tor}^S_1(\Q, R) & \colon \quad 6,\ 8,\ 8,\ 10 , \\
			\operatorname{Tor}^S_2(\Q, R) & \colon \quad 12,\ 14,\ 14,\ 14,\ 16 , \\
			\operatorname{Tor}^S_3(\Q, R) & \colon \quad 18,\ 20 .
		\end{aligned}
	\end{equation*}
\end{proposition}

\begin{proof}
	The complex $ K $ has dimension $ 7 \cdot 8 = 56 $ over $ \Q $, and it is the direct sum of its parts of fixed internal degree.
	So its homology is computed by linear algebra over $ \Q $.
	In each internal degree one writes the differential as matrices in the basis of $ R $ above, and the dimension of the homology is the dimension of the kernel of one matrix minus the rank of the next.
	A direct computation gives the list above.

	The result passes the following check.
	In each internal degree, the alternating sum over $ i $ of the dimensions is the same for $ K $ and for its homology.
	For $ K $ it is the coefficient of $ q^j $ in the product of \mbox{$ (1-q)(1-q^2)(1-q^3) $} with the Hilbert series of $ R $, in which $ q^j $ records the degree $ 2j $.
	For the homology it is read off from the list above.
	The two agree:
	\begin{equation*}
		(1 + q + 2q^2 + 2q^3 + q^4)(1-q)(1-q^2)(1-q^3) = 1 - (q^3 + 2q^4 + q^5) + (q^6 + 3q^7 + q^8) - (q^9 + q^{10}) . \qedhere
	\end{equation*}
\end{proof}

Next, let us compare this homology with $ H^*(\UU) $.
We give a class of homological degree $ i $ and internal degree $ 2j $ the degree $ 2j - i $.
The degrees of the twelve classes of \autoref{prop:koszul_dimensions} are then
\begin{equation*}
	0, \qquad 5,\ 7,\ 7,\ 9, \qquad 10,\ 12,\ 12,\ 12,\ 14, \qquad 15,\ 17 ,
\end{equation*}
and these are the degrees of a basis of $ H^*(\UU) $.

\begin{proposition}
	\label{prop:koszul_comparison}
	The homology $ \operatorname{Tor}^S(\Q, R) $ and the cohomology $ H^*(\UU) $ agree in the following sense.
	\begin{enumerate}[(1)]
		\item For all $ i $ and $ j $,
		\begin{equation*}
			\dim \operatorname{Tor}^S_i(\Q, R)_{2j} = \text{the multiplicity of $ \Q(-j) $ in $ H^{2j-i}(\UU) $} .
		\end{equation*}
		\item If $ \operatorname{Tor}^S_i(\Q, R)_{2j} \neq 0 $, then $ H^{2j-i}(\UU) $ is pure of weight $ 2j $.
		So the homological degree is the weight of a class of $ H^*(\UU) $ minus its degree.
		\item The algebra $ \operatorname{Tor}^S(\Q, R) $ has a basis in which its multiplication table is that of \autoref{thm:eqs_mult}.
		The generators $ x $, $ y $, $ y' $, $ z $ of that table have homological degree $ 1 $ and internal degrees $ 6 $, $ 8 $, $ 8 $, $ 10 $, and the generators $ t $, $ e $ have homological degree $ 2 $ and internal degrees $ 12 $, $ 14 $.
	\end{enumerate}
\end{proposition}

\begin{proof}
	Part (1) is the comparison of \autoref{prop:koszul_dimensions} with the table of \autoref{thm:cohom_of_equations}.
	Part (2) follows from it, because every $ H^n(\UU) $ in that table is a sum of copies of a single $ \Q(-j) $, which has weight $ 2j $.

	Let us prove part (3).
	The products in $ \operatorname{Tor}^S(\Q, R) $ are computed in the same way, on representing cycles.
	Let $ x $ be a generator of $ \operatorname{Tor}^S_1(\Q, R)_{6} $.
	A direct computation shows the following.
	Multiplication by $ x $ is injective on $ \operatorname{Tor}^S_1(\Q, R)_{8} $, on $ \operatorname{Tor}^S_1(\Q, R)_{10} $ and on $ \operatorname{Tor}^S_2(\Q, R)_{12} $, and it maps $ \operatorname{Tor}^S_2(\Q, R)_{14} $ onto $ \operatorname{Tor}^S_3(\Q, R)_{20} $.
	The products of two classes of $ \operatorname{Tor}^S_1(\Q, R)_{8} $, and of a class of $ \operatorname{Tor}^S_1(\Q, R)_{8} $ with a class of $ \operatorname{Tor}^S_2(\Q, R)_{12} $, are zero.

	Choose a basis $ y $, $ y' $ of $ \operatorname{Tor}^S_1(\Q, R)_{8} $, and generators $ z $ of $ \operatorname{Tor}^S_1(\Q, R)_{10} $ and $ t $ of $ \operatorname{Tor}^S_2(\Q, R)_{12} $.
	Then $ xy $ and $ xy' $ are linearly independent in $ \operatorname{Tor}^S_2(\Q, R)_{14} $, which has dimension $ 3 $, and $ xz $ and $ xt $ are non-zero.
	Since $ x^2 = 0 $, the classes $ xy $ and $ xy' $ span the kernel of the multiplication by $ x $ on $ \operatorname{Tor}^S_2(\Q, R)_{14} $.
	Choose a class $ e $ of $ \operatorname{Tor}^S_2(\Q, R)_{14} $ outside this kernel, so that $ xe \neq 0 $, as in \autoref{rem:e_normalisation}.
	The twelve classes
	\begin{equation*}
		1, \quad x, \quad y, \quad y', \quad z, \quad t, \quad e, \quad xy, \quad xy', \quad xz, \quad xt, \quad xe
	\end{equation*}
	form a basis of $ \operatorname{Tor}^S(\Q, R) $ by \autoref{prop:koszul_dimensions}.
	Every product of two of them that is not listed in the table of \autoref{thm:eqs_mult} is zero.
	Indeed, it either contains $ x^2 $, or is one of the two products found above to be zero, or lies in a part of $ \operatorname{Tor}^S(\Q, R) $ that is zero by \autoref{prop:koszul_dimensions}.
	The signs are those of the table, because two classes of $ \operatorname{Tor}^S(\Q, R) $ anticommute exactly when both have odd homological degree, that is, odd degree.
\end{proof}

We observe part (2) by comparing two computations and do not derive it from mixed Hodge theory.

The two sides of this comparison are related by a spectral sequence.

\begin{proposition}
	\label{prop:koszul_eilenberg_moore}
	The Eilenberg--Moore spectral sequence of the fibre bundle $ \UU \to \UU_{hG} \to BG $ degenerates at the second page.
\end{proposition}

\begin{proof}
	The space $ \UU $ is the fibre of this bundle, and $ H^*(\UU_{hG}) = R $.
	So the spectral sequence converges to $ H^*(\UU) $, and its second page is $ \operatorname{Tor}^S(\Q, R) $, with $ \operatorname{Tor}^S_i(\Q, R)_{2j} $ in total degree $ 2j - i $.
	Its limit page is a subquotient of its second page.
	Both $ \operatorname{Tor}^S(\Q, R) $ and $ H^*(\UU) $ have total dimension $ 12 $, by \autoref{prop:koszul_dimensions} and \autoref{thm:cohom_of_equations}.
	So the second page is the limit page.
\end{proof}

Finally, let us say what the generators of \autoref{thm:eqs_mult} mean in commutative algebra.
The socle of $ R $ is the space of elements of $ R $ annihilated by $ s_1 $, $ s_2 $ and $ s_3 $.

\begin{proposition}
	\label{prop:koszul_commutative_algebra}
	The ideal $ I $ and the ring $ R = S/I $ have the following properties.
	\begin{enumerate}[(1)]
		\item The ideal $ I $ has four minimal generators, of degrees $ 6 $, $ 8 $, $ 8 $ and $ 10 $.
		\item It has five minimal first syzygies, of degrees $ 12 $, $ 14 $, $ 14 $, $ 14 $ and $ 16 $.
		\item The socle of $ R $ has dimension $ 2 $, and it is spanned by
		\begin{equation*}
			51 s_1^3 - 25 s_1 s_2 \qquad \text{and} \qquad s_2^2 .
		\end{equation*}
		\item The ring $ R $ is neither a complete intersection nor a Gorenstein ring.
	\end{enumerate}
\end{proposition}

\begin{proof}
	We use \autoref{lem:koszul_tor} and \autoref{prop:koszul_dimensions}.

	The space $ \operatorname{Tor}^S_1(\Q, R) $ is $ I / (s_1,s_2,s_3) I $, the space of minimal generators of $ I $.
	This gives part (1).

	The space $ \operatorname{Tor}^S_2(\Q, R) $ is $ \operatorname{Tor}^S_1(\Q, I) $, the space of minimal first syzygies of $ I $, that is, of the minimal relations between its minimal generators.
	This gives part (2).

	The space $ \operatorname{Tor}^S_3(\Q, R) $ is the socle of $ R $, with the internal degree raised by $ 12 $.
	Indeed, an element $ r \sigma_1 \sigma_2 \sigma_3 $ of $ K $ is a cycle exactly when $ s_1 r = s_2 r = s_3 r = 0 $, and no such element is a boundary.
	So the socle has dimension $ 2 $, with one generator in degree $ 6 $ and one in degree $ 8 $.
	The element $ s_2^2 $ lies in it, because $ R $ is zero in degrees above $ 8 $.
	So does $ 51 s_1^3 - 25 s_1 s_2 $.
	Indeed, its products with $ s_2 $ and $ s_3 $ have degrees above $ 8 $, and by the relations of degree $ 8 $
	\begin{equation*}
		s_1 (51 s_1^3 - 25 s_1 s_2) = \Bigl( 51 \cdot \frac{25}{369} - 25 \cdot \frac{17}{123} \Bigr) s_2^2 = 0 .
	\end{equation*}
	This gives part (3).

	The ideal of a complete intersection would have three minimal generators, and the socle of a Gorenstein ring would have dimension $ 1 $.
	This gives part (4).
\end{proof}

By part (1), one of the three generators of degree $ 8 $ in \autoref{cor:cohom_ring_moduli_presentation} is redundant.
Indeed, the product of $ s_1 $ with the generator of degree $ 6 $ is the linear combination of these three with the coefficients $ 72 $, $ -39 $ and $ 7 $.

In the basis of \autoref{prop:koszul_comparison}, the classes $ x $, $ y $, $ y' $, $ z $ correspond to the minimal generators of $ I $.
The classes $ t $, $ e $, $ xy $, $ xy' $, $ xz $ correspond to its minimal first syzygies, and the classes $ xt $, $ xe $ to the socle of $ R $.

\printbibliography

\end{document}